\documentclass[mreferee]{gji}
\usepackage{footnote}
\usepackage{timet}
\usepackage{algorithm}
\usepackage{algorithmic}
\usepackage[pdftex]{graphicx}
\usepackage{epstopdf}
\usepackage{epsfig,subfigure}
\usepackage{epsf}
\usepackage{graphics}
\usepackage{url}
\usepackage{color}
\usepackage{amssymb}
\usepackage{amsmath}

\graphicspath{{./Figures/}}

\addtocounter{footnote}{0}
\newcommand{\bfm}{\mathbf{m}}
\newcommand{\bfd}{\mathbf{d}}

\newcommand{\bfr}{\mathbf{r}}
\newcommand{\bfy}{\mathbf{y}}
\newcommand{\bfh}{\mathbf{h}}

\newcommand{\bfu}{\mathbf{u}}

\newcommand{\bfdo}{\mathbf{d}_{\mathrm{obs}}}
\newcommand{\bfdp}{\mathbf{d}_{\mathrm{pre}}}

\newcommand{\bfde}{\mathbf{d}_{\mathrm{exact}}}
\newcommand{\bfma}{\mathbf{m}_{\mathrm{apr}}}

\newcommand{\Wd}{W_{\bfd}}

\newcommand{\Wde}{W_{\mathrm{depth}}}
\newcommand{\WR}{W_{\mathrm{R}}}
\newcommand{\WRx}{W_{\mathrm{R_x}}}
\newcommand{\WRy}{W_{\mathrm{R_y}}}
\newcommand{\WRz}{W_{\mathrm{R_z}}}
\newcommand{\Wr}{w_{\mathrm{r}}}
\newcommand{\Dx}{D_{\mathrm{x}}}
\newcommand{\Dy}{D_{\mathrm{y}}}
\newcommand{\Dz}{D_{\mathrm{z}}}
\newcommand{\px}{p_{\mathrm{x}}}
\newcommand{\py}{p_{\mathrm{y}}}
\newcommand{\pz}{p_{\mathrm{z}}}

\newcommand{\diag}{\mathrm{diag}}

\newcommand{\Rn}{\mathcal{R}^{n}}
\newcommand{\Rm}{\mathcal{R}^{m}}
\newcommand{\Rmn}{\mathcal{R}^{m \times n}}
\newcommand{\Rmq}{\mathcal{R}^{m \times q}}
\newcommand{\Rqq}{\mathcal{R}^{q \times q}}

\title[TV regularization of the gravity inverse problem using RGSVD]{Total variation regularization of the $3$-D gravity inverse problem using a randomized generalized singular value decomposition}
\author[S. Vatankhah, R.~A. Renaut, V.~E. Ardestani]{Saeed Vatankhah $^1$, Rosemary A. Renaut $^2$ and Vahid E. Ardestani $^1$ \\
$^1$ Institute of Geophysics, University of Tehran, Iran \\
$^2$ School of Mathematical and Statistical Sciences, Arizona State University, Tempe, AZ, USA.}
\begin{document}
\maketitle

\begin{summary}
We present a fast algorithm for the total variation regularization of the $3$-D gravity inverse problem. Through imposition of the total variation regularization, subsurface structures presenting with sharp discontinuities  are preserved better than when using a conventional minimum-structure inversion. The associated problem formulation for the regularization is non linear but can be solved using an iteratively reweighted least squares  algorithm. For small scale problems the regularized least squares problem at each iteration can be solved using the generalized singular value decomposition. This is not feasible for large scale problems. Instead we introduce the use of a randomized generalized singular value decomposition in order to reduce the dimensions of the problem and provide an effective and efficient solution technique. For further efficiency an alternating direction algorithm is used to implement the total variation weighting operator within the iteratively reweighted least squares algorithm.   Presented results for synthetic examples demonstrate that the novel randomized decomposition provides good accuracy for reduced computational and memory demands as compared to use of classical approaches. 
  
\end{summary}

\section{Introduction}\label{sec:intro}
Regularization is imposed in order to find acceptable solutions to the ill-posed and non-unique problem of gravity inversion.  Most current regularization techniques minimize  a global objective function that consists of a data misfit term and a stabilization term; \cite{LiOl:98,PoZh:99,BoCh:2001,VRA:2017a}. Generally for potential field inversion the data misfit  is measured as a weighted $L_2$-norm of the difference between the observed and predicted data of the reconstructed model, \cite{Pi:09}. Stabilization aims to both counteract the ill-posedness of the problem so that changes in the model parameters due to small changes in the observed data are controlled, and to impose realistic characteristics  on the reconstructed model. Many forms of robust and reliable stabilizations have been used by the geophysical community.  $L_0$, $L_1$ and Cauchy norms for the model parameters yield sparse and compact solutions, \cite{LaKu:83,PoZh:99,Ajo:2007,Pi:09,VRA:2017a}.  The minimum-structure inversion based on a $L_2$ measure of the gradient of the model parameters  yields a model that is smooth \cite{LiOl:98}. Total variation (TV) regularization based on a  $L_1$ norm of the gradient of the model parameters \cite{BCO:2002,Far:2008}  preserves edges in the model and provides a reconstruction of the model that is blocky and non-smooth. While the selection of stabilization for a given data sets depends on an anticipated true representation of the subsurface model and is problem specific,  efficient and practical algorithms for large scale problems are desired for all formulations. Here the focus is on the development of a new randomized algorithm for the $3$-D linear inversion of gravity data with TV stabilization. 

The objective function to be minimized through the stabilization techniques used for geophysical inversion is nonlinear in the model and a solution is typically found using an iteratively reweighted standard Tikhonov least squares (IRLS) formulation. For the TV constraint  it is necessary to solve a generalized Tikhonov least squares problem at each iteration. This presents no difficulty for small scale problems; the generalized singular value decomposition (GSVD) can be used to compute the mutual decomposition of the model and  stabilization matrices from which a solution is immediate. Furthermore, the use of the GSVD provides a convenient form for the regularization parameter-choice rules,  \cite{XiZo:2013,ChPa:2015}. The approach using the GSVD is nor practical for large scale problems, neither in terms of computational costs nor memory demands.  Randomized algorithms  compute low-rank matrix approximations with reduced memory and computational costs,  \cite{Halko:2011,XiZo:2013,VMN:2015,WXZ:2016}. Random sampling is used to construct a low-dimensional subspace that captures the dominant spectral properties of the original matrix  and then,  as developed by Wei et al. \shortcite{WXZ:2016}, a GSVD can be applied to regularize for this dominant spectral space. We demonstrate that applying this randomized GSVD (RGSVD) methodology at each step of the iterative TV regularization provides a fast and efficient algorithm for $3$-D gravity inversion that inherits all the properties of TV inversion for small-scale problems. We also show how the weighting matrix at each iteration can be determined and how an alternating direction algorithm  can also reduce the memory overhead associated with calculating the GSVD of the projected subproblem at each step of the IRLS. 

\section{ Inversion methodology }\label{inversionmethod}

Suppose the subsurface is divided into a large number of cells of fixed size but unknown density \cite{LiOl:98,BoCh:2001}. Unknown densities of the cells are stacked in vector $\bfm \in \Rn$ and measured data on the surface are stacked in $\bfdo \in \Rm$ are related to the densities via the linear relationship
\begin{eqnarray}\label{d=gm}
\bfdo= G\bfm,
\end{eqnarray}
for forward model matrix $G \in \Rmn$ ($m\ll n$). The aim is to find an acceptable model for the densities that predicts the observed data at the noise level. For gravity inversion,  \eqref{d=gm} is modified through inclusion of an estimated prior model, $\bfma$, and a depth weighting matrix, $\Wde$, \cite{LiOl:98,BoCh:2001} yielding  
\begin{eqnarray}\label{modified d=gm}
\bfdo - G\bfma = G\bfm - G\bfma. 
\end{eqnarray}    
This is replaced first by $G\bfy=\bfr$ using $\bfr=\bfdo - G\bfma$ and $ \bfy= \bfm - \bfma$, and then by incorporating depth weighting,  via $\tilde{G}=G\Wde^{-1}$ and $\bfh=\Wde\bfy$, further replaced by 
\begin{eqnarray}\label{r=gh}
\bfr= \tilde{G}\bfh. 
\end{eqnarray}
Due to the ill-posedness of the problem  \eqref{r=gh} cannot be solved directly but must be stabilized in order to provide a practically acceptable solution. Using the $L_2$-TV regularization methodology $\bfh$ is obtained by minimizing the objective function \cite{WoRo:07},
\begin{eqnarray}\label{globalfunction1}
P^{\alpha}(\bfh)=\| \Wd(\tilde{G}\bfh-\bfr)  \|_2^2 + \alpha^2 \| \vert \nabla \bfh \vert\|_1^1,
\end{eqnarray} 
in which  diagonal weighting matrix $\Wd$ approximates the square root of the inverse covariance matrix for the independent noise in the data. Specifically $(\Wd)_{ii}= 1/\eta_i$  in which $\eta_i $ is the standard deviation of the noise for the $i$th datum.  Regularization parameter $\alpha$ provides a tradeoff between the weighted data misfit and the stabilization.  The absolute value of the gradient vector, $ \vert \nabla \bfh \vert $, is equal to $= \sqrt{(\Dx\bfh)^2+(\Dy \bfh)^2+(\Dz \bfh)^2}$ in which $\Dx$, $\Dy$ and $\Dz$ are the discrete derivative operators in $x$, $y$ and $z$-directions, respectively. The derivatives at the centers of the cells are approximated to low order using forward differences with backward differencing at the boundary points,  Li $\&$ Oldenburg \shortcite{LiOl:2000},  yielding matrices $\Dx$, $\Dy$ and $\Dz$  which are square of size $n \times n$.  Although \eqref{globalfunction1} is a convex optimization problem with a unique solution, the TV term is not differentiable everywhere with respect to $\bfh$ and to find a solution it is helpful to rewrite the stabilization term using a weighted $L_2$-norm, following Wohlberg \& Rodr{\'i}guez \shortcite{WoRo:07}. Given vectors $\xi$, $\chi$ and $\psi$ and a diagonal matrix $\WR$ with entries $\Wr$ we have
\begin{eqnarray}\label{L2normVectors}
\biggr\Vert \left(  \begin{tabular}{c c c }  $\WR$ & $0$ & $0$ \\  $0$ & $\WR$ & $0$  \\  $0$ & $0$ & $\WR$    \end{tabular} \right)  \left(  \begin{tabular}{c }  $\xi$ \\ $\chi$ \\ $\psi$      \end{tabular} \right)   \biggr \Vert_2^2 =
\sum_r(\Wr^2 \xi_{\mathrm{r}}^2 + \Wr^2 \chi_{\mathrm{r}}^2 + \Wr^2 \psi_{\mathrm{r}}^2) =
\sum_r \Wr^2 (\sqrt{\xi_{\mathrm{r}}^2 + \chi_{\mathrm{r}}^2 + \psi_{\mathrm{r}}^2})^2.
\end{eqnarray}
Setting $w_{\mathrm{r}}= (\xi_{\mathrm{r}}^2 + \chi_{\mathrm{r}}^2 + \psi_{\mathrm{r}}^2)^{-1/4} $ we have
\begin{eqnarray}\label{L2toL1}
\biggr \Vert \left(  \begin{tabular}{c c c }  $\WR$ & $0$ & $0$ \\  $0$ & $\WR$ & $0$  \\  $0$ & $0$ & $\WR$    \end{tabular} \right)  \left(  \begin{tabular}{c }  $\xi$ \\ $\chi$ \\ $\psi$      \end{tabular} \right)   \biggr \Vert_2^2 =\sum_r (\sqrt{\xi_{\mathrm{r}}^2 + \chi_{\mathrm{r}}^2 + \psi_{\mathrm{r}}^2}).
\end{eqnarray}
Now with $\xi=\Dx\bfh$, $\chi=\Dy\bfh$, $\psi=\Dz\bfh$ and  $\WR^{(k)}$ defined to have entries $\Wr^{(k)}$ calculated for $\nabla \bfh$ at iteration $k-1$
given by $\Wr^{(k)}= ((\Dx\bfh^{(k-1)})_r^2 + (\Dy\bfh^{(k-1)})_r^2 + (\Dz\bfh^{(k-1)})_r^2 +\epsilon^2)^{-1/4}$, the TV stabilizer is approximated via
\begin{eqnarray}
\Vert \vert \nabla \bfh \vert \Vert_1^1= \Vert \sqrt{(\Dx\bfh)^2+(\Dy \bfh)^2+(\Dz \bfh)^2} \Vert_1^1 \approx \Vert W^{(k)} D \bfh \Vert_2^2,
\end{eqnarray} 
for derivative operator $D=[\Dx; \Dy;\Dz]$. 
Here, $0<\epsilon\ll1$  is added to avoid the possibility of division by zero, and 
superscript $k$ indicates that  matrix $\WR$ is updated using the model parameters of the previous iteration. Hence \eqref{globalfunction1} is rewritten as a general Tikhonov functional 
\begin{eqnarray}\label{globalfunction2}
P^{\alpha}(\bfh)=\| \Wd(\tilde{G}\bfh-\bfr)  \|_2^2 + \alpha^2 \| \tilde{D} \bfh \|_2^2, \quad \tilde{D}=WD,
\end{eqnarray} 
for which the minimum is explicitly expressible as 
\begin{eqnarray}\label{hsolution}
\bfh=(\tilde{\tilde{G}}^T\tilde{\tilde{G}}+\alpha^2\tilde{D}^T\tilde{D})^{-1}\tilde{\tilde{G}}^T\tilde{\bfr}, \quad \tilde{\tilde{G}}=\Wd\tilde{G} \quad  \mathrm{and} \quad \tilde{\bfr}=\Wd\bfr.
\end{eqnarray}
The model update is then
\begin{eqnarray}\label{modelupdate}
\bfm(\alpha)=\bfma+\Wde^{-1}\bfh (\alpha)
\end{eqnarray}
While the Tikhonov function can be replaced by a standard form Tikhonov function, i.e. with regularization term $\|\bfh\|^2$, when $\tilde{D}$ is easily invertible, e.g. when diagonal,  the form of $\tilde{D}$ in this case makes that transformation prohibitive for cost and we must solve using the general form.

We show the iterative process for the solution in Algorithm~\ref{IterativeTVRGSVD}. It should be noted  that the iteration process terminates when solution satisfies the noise level or a predefined maximum number of iterations is reached \cite{BoCh:2001}.  Furthermore, the positivity constraint [$\rho_{\mathrm{min}}, \rho_{\mathrm{max}}$] is imposed at each iteration. If at any iteration a density value falls outside these predefined density bounds, the value is projected back to the nearest bound value \cite{BoCh:2001}. 

For small-scale problems in which the dimensions of $\tilde{\tilde{G}}$, and consequently $\tilde{D}$, are small, the the solution $\bfh(\alpha)$ is found at minimal cost using the GSVD of the matrix pair $[\tilde{\tilde{G}},\tilde{D}]$ as it is shown in  Appendix~\ref{gsvd} \cite{ABT:2013,VAR:2014}. Furthermore, given the GSVD the regularization parameter may be estimated cheaply using standard parameter-choice techniques \cite{XiZo:2013,ChPa:2015}. But for large scale problems it is not practical to calculate the GSVD at each iteration, both with respect to computational cost and memory demands. Instead the size of the original large problem can be reduced greatly using a randomization technique which provides the  GSVD in a more feasible and efficient manner.  The solution of reduced system still is a good approximation of the original system  \cite{Halko:2011,XiZo:2013,VMN:2015,XiZo:2015,WXZ:2016}. Here, we use the Randomized GSVD (RGSVD) algorithm developed by Wei et al. \shortcite{WXZ:2016} for under-determined problems, in which the mutual decomposition of the matrix pair $[\tilde{\tilde{G}},\tilde{D}]$ is approximated by 
\begin{equation}\label{GSVD}
\tilde{\tilde{G}} \approx U \Lambda Z, \quad \tilde{D}=V M Z, \quad U \in \Rmq, V \in \mathcal{R}^{3n \times q}, \Lambda \in \Rqq, M \in \Rqq, Z \in \mathcal{R}^{q \times n}.
\end{equation}

\begin{algorithm}
\caption{RGSVD algorithm. Given matrices $\tilde{\tilde{G}} \in \Rmn (m \leq n)$ and $\tilde{D} \in \mathcal{R}^{3n \times n}$, a target matrix rank $q $ and a small constant oversampling parameter $p$ satisfying $q+p=l \ll m$, compute an approximate GSVD of $[\tilde{\tilde{G}},\tilde{D}]$: $\tilde{\tilde{G}} \approx U \Lambda Z$, $\tilde{D}=V M Z$ with $U \in \Rmq$, $V \in \mathcal{R}^{3n \times q}$, $\Lambda \in \Rqq$, $M \in \Rqq$ and $Z \in \mathcal{R}^{q \times n}$.}\label{RGSVDAlgorithm}
\begin{algorithmic}[1]
\STATE Generate a Gaussian random matrix $\Omega \in \mathcal{R}^{l \times m} $.
\STATE Form the matrix $Y=\Omega \tilde{\tilde{G}} \in \mathcal{R}^{l \times n}$.  
\STATE Compute orthonormal matrix $Q \in \mathcal{R}^{n \times l}$ via QR factorization $Y^T=QR$.
\STATE Set $Q=Q(:,1:q)$ and form the matrices $B_1=\tilde{\tilde{G}}Q \in \mathcal{R}^{m \times q}$ and $B_2=\tilde{D}Q \in \mathcal{R}^{3n \times q}$.
\STATE Compute the GSVD of $[B_1, B_2]$: $ \left[ \begin{tabular}{c} $B_1$ \\ $B_2$ \end{tabular} \right]$ $=$ $ \left[ \begin{tabular}{c c} $U$ & \\  & $V$ \end{tabular} \right]$ $ \left[ \begin{tabular}{c} $\Lambda$ \\ $M$ \end{tabular} \right]$ $X^T$, using $[U,V,X,\Lambda,M]=\text{gsvd}(B_1,B_2,0)$. 
\STATE Form the matrix  $Z=X^T Q^T \in \mathcal{R}^{q \times n}$.
\end{algorithmic}
\end{algorithm}

The steps of the algorithm are given in Algorithm~\ref{RGSVDAlgorithm}. Steps $1$ to $3$ are used to form matrix $Q$ which approximates the range of $\tilde{\tilde{G}}^T$. At step $4$,  $\tilde{\tilde{G}}$ is projected into a lower dimensional matrix $B_1$, for which $B_1$ provides information on the range of $\tilde{\tilde{G}}$. The same projection is applied to matrix $\tilde{D}$. In step $5$, an economy-sized GSVD is computed for the matrix pair $[B_1,B_2]$.  Parameter $q$ balances the accuracy and efficiency of the Algorithm~\ref{RGSVDAlgorithm} and determines the dimension of the subspace for the projected problem. When $q$ is small, this methodology is very effective and leads to a fast GSVD computation. Simultaneously, the parameter $q$ has to be selected large enough to capture the dominant spectral properties of the original problem with the aim that the  solution obtained using the RGSVD is close to the solution that would be obtained using the GSVD.

For any GSVD decomposition \eqref{GSVD}, including that obtained via Algorithm~\ref{RGSVDAlgorithm}, $\bfh(\alpha)$ of \eqref{hsolution} is given by
\begin{eqnarray}\label{rgsvdsoln1}
\bfh(\alpha)= \left( Z^T \Lambda^T U^T U \Lambda Z + \alpha^2 Z^T M^T V^T V M Z \right)^{-1} Z^T \Lambda^T U^T \tilde{\bfr},
\end{eqnarray}
which simplifies to, \cite{ABT:2013,VAR:2014}, 
\begin{eqnarray}\label{rgsvdsoln2}
\bfh(\alpha)= \sum_{i=1}^{q} \frac{\gamma^2_i}{\gamma^2_i+\alpha^2} \frac{\bfu^T_{i}\tilde{\bfr}}{\lambda_i} (Z^{-1})_{i}.
\end{eqnarray}
Here $\gamma_i$ is the $i$th generalized singular value, see Appendix~\ref{gsvd}, and $Z^{-1}$ is the Moore-Penrose inverse of $Z$ \cite{WXZ:2016}. Incorporating this single step of the RGSVD within the TV algorithm, yields the iteratively reweighted TVRGSVD algorithm given in Algorithm~\ref{IterativeTVRGSVD}. Note here that the  steps $1$-$3$ and the calculation of $B_1$ in step $4$ in Algorithm~\ref{RGSVDAlgorithm} are the same for all TV iterations and thus outside the loop in Algorithm~\ref{IterativeTVRGSVD}. At each iteration $k$, the matrices $W$, $\tilde{D}$ and $B_2$ are updated, and the  GSVD is determined for the matrix pair $[B_1, B_2^{(k)}]$. We should note here that it is immediate to use the  minimum gradient support (MGS) stabilizer introduced by Portniaguine \& Zhdanov \shortcite{PoZh:99} in Algorithm~\ref{IterativeTVRGSVD} via replacing $(-1/4)$ in $\WR$ with $(-1/2)$ and keeping all other parameters are kept fixed.

\begin{algorithm}
\caption{Iterative TV inversion algorithm using randomized GSVD}\label{IterativeTVRGSVD}
\begin{algorithmic}[1]
\REQUIRE $\bfdo$, $\bfma$, $G$, $\Wd$, $\Wde$, $\Dx$, $\Dy$, $\Dz$, $q$, $\epsilon > 0$, $\rho_{\mathrm{min}}$, $\rho_{\mathrm{max}}$, $K_{\mathrm{max}}$
\STATE Initialize $\bfm^{(0)}=\bfma$, $W^{(1)} =I$, $D=[\Dx;\Dy;\Dz]$, $\tilde{D}^{(1)}=D$, $k=1$
\STATE Calculate $\tilde{\tilde{G}}=\Wd G \Wde^{-1}$, $\tilde{\bfr}^{(1)}=\Wd(\bfdo-G\bfma)$
\STATE Generate a Gaussian random matrix $\Omega \in \mathcal{R}^{l \times m} $.
\STATE Compute matrix $Y=\Omega \tilde{\tilde{G}} \in \mathcal{R}^{l \times n}$.  
\STATE Compute orthonormal matrix $Q \in \mathcal{R}^{n \times l}$ via QR factorization $Y^T=QR$.
\STATE Set $Q=Q(:,1:q)$ and form the the matrix $B_{1}=\tilde{\tilde{G}}Q \in \mathcal{R}^{m \times q}$.
\WHILE {Not converged, noise level not satisfied, and $k<K_{\mathrm{max}}$} 
\STATE Form the matrix $B_{2}=\tilde{D}^{(k)}Q \in \mathcal{R}^{3n \times q}$.
\STATE \label{gsvdstep} Compute the GSVD of $\lbrace B_1,B_2 \rbrace$:$[U,V,X,\Lambda,M]=\text{gsvd}(B_1,B_2,0)$.
\STATE Form the matrix $Z=X^TQ^T \in \mathcal{R}^{q \times n}$. 
\STATE \label{stepalpha} Estimate $\alpha^{(k)}$ using \eqref{upregsvd} 
\STATE Set $\bfh^{(k)}= \sum_{i=1}^{q} \frac{\gamma_i^2}{\gamma_i^2+(\alpha^{(k)})^2} \frac{\bfu_i ^T\tilde{\bfr}^{(k)}}{\lambda_i} (Z^{-1})_i$.
\STATE Set $\bfm^{(k)}=\bfm^{(k-1)}+ (\Wde)^{-1}\bfh^{(k)}$.
\STATE Impose constraint conditions on $\bfm^{(k)}$ to force $\rho_{\mathrm{min}}\le \bfm^{(k)} \le \rho_{\mathrm{max}}$.
\STATE Test convergence and exit loop if converged.
\STATE \label{Wupdate} Calculate $\WR^{(k+1)}=\diag\left(  ((\Dx\bfh^{(k)})^2 + (\Dy\bfh^{(k)})^2 + (\Dz\bfh^{(k)})^2 +\epsilon^2)^{-1/4}\right)$ and set $W^{(k+1)}=\diag(\WR^{(k+1)};\WR^{(k+1)};\WR^{(k+1)})$.
\STATE \label{Dupdate} Set  $\tilde{\bfr}^{(k+1)}=\Wd(\bfdo-G\bfm^{(k)})$, $\tilde{D}^{(k+1)}=W^{(k+1)}D$, 
\STATE {$k=k+1$}
\ENDWHILE
\ENSURE Solution $\rho=\bfm^{(k)}$. $K=k$.
\end{algorithmic}
\end{algorithm}

\subsubsection{An Alternating Direction Algorithm}
Step \ref{gsvdstep}  of Algorithm~\ref{IterativeTVRGSVD} requires the economy GSVD decomposition in which matrix $B_2^{(k)}$ is of size $3n \times q$. For a large scale problem the computational cost and memory demands with the calculation of the GSVD limits the size of the problem that can be solved. We therefore turn to an alternative approach for large scale three dimensional problems and adopt the use of an Alternating Direction (AD) strategy, in which to handle the large scale problem requiring derivatives in greater than two dimensions we split the problem into pieces handling each direction one after the other. This is a technique that has been in the literature for some time for handling the solution of large scale partial differential equations, most notably through the alternating direction implicit method \cite{ADM}. Matrix  $D^{(k)}$ generated  via steps~\ref{Wupdate} and \ref{Dupdate} of Algorithm~\ref{IterativeTVRGSVD}  can be changed without changing the other steps of the algorithm. For the AD algorithm we may therefore alternate over $D^{(k)} = \WRx^{(k)}\Dx$, $\WRy^{(k)}\Dy$ or $\WRz^{(k)}\Dz$, dependent on $(k\mod 3)=0$, $1$, or $2$, respectively. Then $B_2^{(k)}$ is only of size $n \times n$, yielding reduced memory demands for calculating the GSVD. We note that $D$ is also initialized consistently. The AD approach amounts to apply the edge preserving algorithm in each direction independently and cycling over all directions. Practically, we also find that there is nothing to be gained by requiring that the derivative matrices are square, and following \cite{Hansen:2007} ignore the derivatives at the boundaries. For a one dimensional problem this amounts to taking matrix $\Dx$ of size $(n-1) \times n$ for a line with $n$ points. Thus dependent on $k \mod 3$ we have a weighting matrix of size $\px \times \px$, $\py \times \py$ or $\pz \times \pz$ for matrices $\Dx$, $\Dy$ and $\Dz$ of sizes $\px \times n$, $\py \times n$ and $\pz \times n$, respectively, and where $\px$, $\py$ and $\pz$ (all less than $n$) depend on the number of boundary points in each direction. Matrix $B_2^{(k)}$ is thus reduced in size and the GSVD calculation is more efficient at each step. Furthermore, we find that rather than calculating the relevant weight matrix $\WR$ for the given dimension we actually form the weighted entry that leads to approximation of the relevant component of the  gradient vector for all three dimensions, thus realistically approximating the gradient as given in \eqref{L2normVectors}. For the presented results we will use the AD version of Algorithm~\ref{IterativeTVRGSVD}, noting that this is not necessary for the smaller problems, but is generally more efficient and reliable for the larger three dimensional problem formulations.

\subsection{Estimation of the Regularization Parameter $\alpha$}
As presented to this point we have assumed a known value for the regularization parameter $\alpha$. Practically we wish to find $\alpha$ dynamically to appropriately regularize at each step of the iteration so as to recognize that the conditioning of the problem changes with the iteration. Here we use the method of unbiased predictive risk estimation (UPRE) which we have found to be robust in our earlier work \cite{RVA:2017,VAR:2015,VRA:2017a}. The method, which goes back to \citename{Mallows} \shortcite{Mallow}, requires some knowledge of the noise level in the data and was carefully developed in Vogel \shortcite{Vogel:2002} for the standard Tikhonov functional.  The method was further extended for use with TV regularization by Lin et al. \shortcite{LWG:2010}. Defining the residual $R(\bfh(\alpha)) = \tilde{\tilde{G}} \bfh(\alpha)-\tilde{\bfr}$ and influence matrix $H_{TV,\alpha}=\tilde{\tilde{G}}(\tilde{\tilde{G}}^T\tilde{\tilde{G}}+\alpha^2\tilde{D}^T\tilde{D})^{-1}\tilde{\tilde{G}}^T$, the optimal parameter $\alpha$ is the minimizer of
\begin{eqnarray}\label{upre}
U(\alpha)=\|R(\bfh(\alpha))\|_2^2+ 2 \,\text{trace}(H_{TV,\alpha})-m,
\end{eqnarray} 
which is given in terms of the GSVD by
\begin{eqnarray}\label{upregsvd}
U(\alpha)=\sum_{i=1}^{q}  \left( \frac{1}{\gamma_i^2 \alpha^{-2} + 1} \right)^2 \left(\bfu_i^T\tilde{\bfr} \right)^2 + 2 \left( \sum_{i=1}^{q} \frac{\gamma_i^2}{\gamma_i^2+\alpha^2}\right) - q.
\end{eqnarray}
Typically $\alpha_{opt}$ is found by evaluating \eqref{upregsvd} on a range of $\alpha$, between minimum and maximum $\gamma_i$, and then that $\alpha$ which minimizes the function is selected as $\alpha_{opt}$.

\section{Synthetic examples}\label{synthetic}

\subsection{Model consisting of two dipping dikes}\label{twodikes}

\begin{figure*}
\subfigure{\label{fig1a}\includegraphics[width=.45\textwidth]{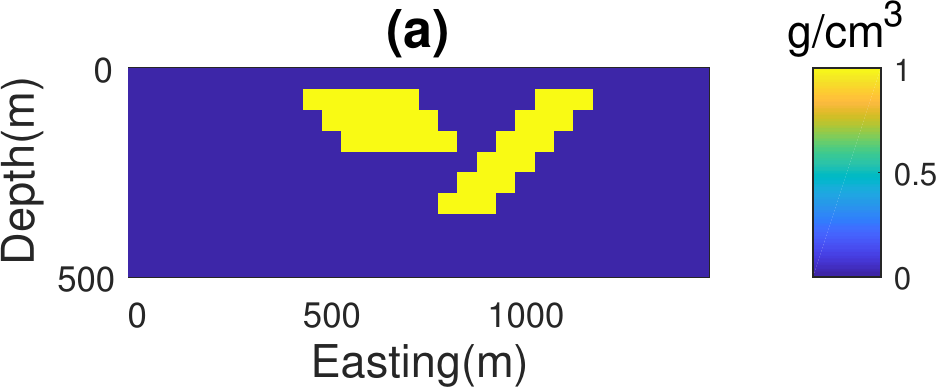}}
\subfigure{\label{fig1b}\includegraphics[width=.35\textwidth]{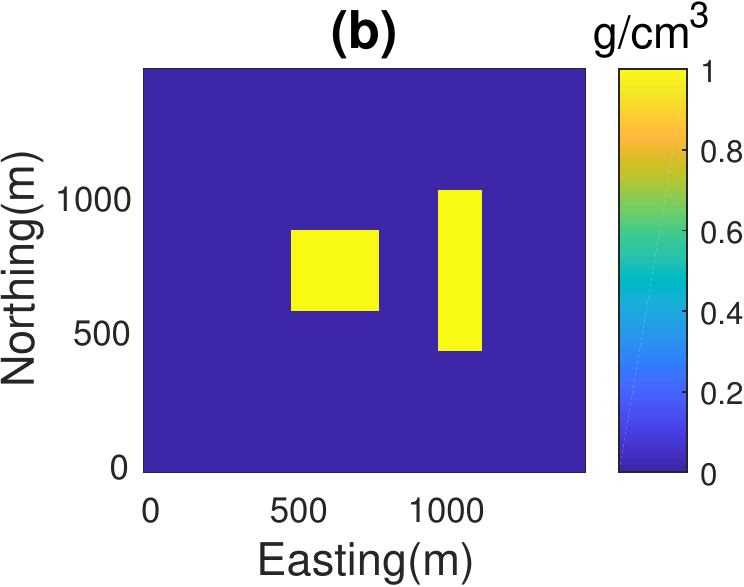}}
\caption {A model that consists of two dikes dipping in opposite directions. (a) Cross-section at  $\mathrm{northing}=725$~m; (b) Plane-section at $\mathrm{depth}=100$~m.} \label{fig1}
\end{figure*}

\begin{figure*}
\includegraphics[width=.35\textwidth]{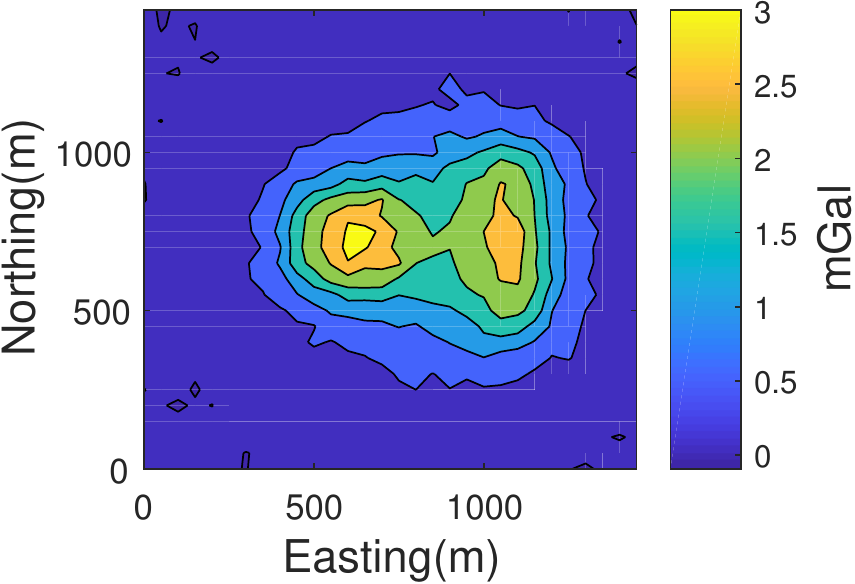}
\caption {The gravity anomaly produced by the model shown in Fig.~\ref{fig1} and contaminated by Gaussian noise.} \label{fig2}
\end{figure*}

As a first example, we use a complex model that consists of two embedded dipping dikes of different sizes and dipping in opposite directions but with the same density contrast $1$~g~cm$^{-3}$, Fig.~\ref{fig1}.  Gravity data, $\bfde$, is generated on the surface for a grid of  $30 \times 30 = 900$ points  with grid spacing $50$~m. Gaussian noise with standard deviation $ (0.02~(\bfde)_i + 0.002~\| \bfde \|)$ is added to each datum. Example noisy data, $\bfdo$, is illustrated in Fig.~\ref{fig2}. Inversion is performed for the subsurface volume of $9000$  cubes of size $50$~m in each dimension using the matrix $\tilde{\tilde{G}}$ of size $900 \times 9000$. Use of this relatively small model permits examination of the inversion methodology with respect to different parameter choices and provides the framework to be used for more realistic larger models. All  computations are performed on a desktop computer with Intel Core i7-4790 CPU 3.6GHz processor and 16 GB RAM. 

\begin{figure*}
\subfigure{\label{fig3a}\includegraphics[width=.45\textwidth]{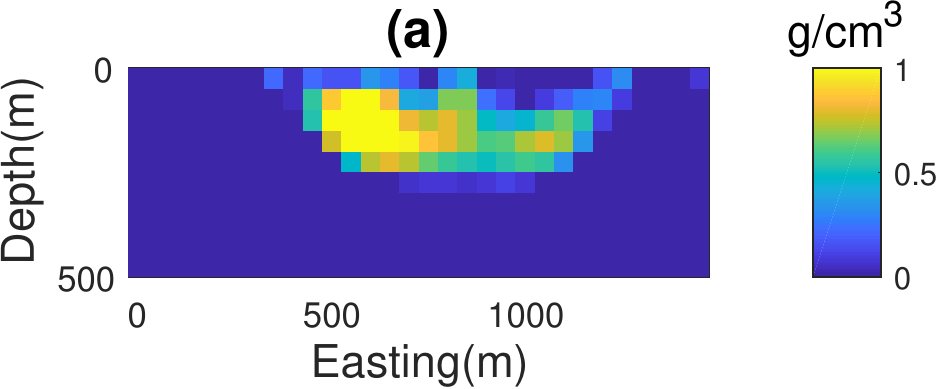}}
\subfigure{\label{fig3b}\includegraphics[width=.35\textwidth]{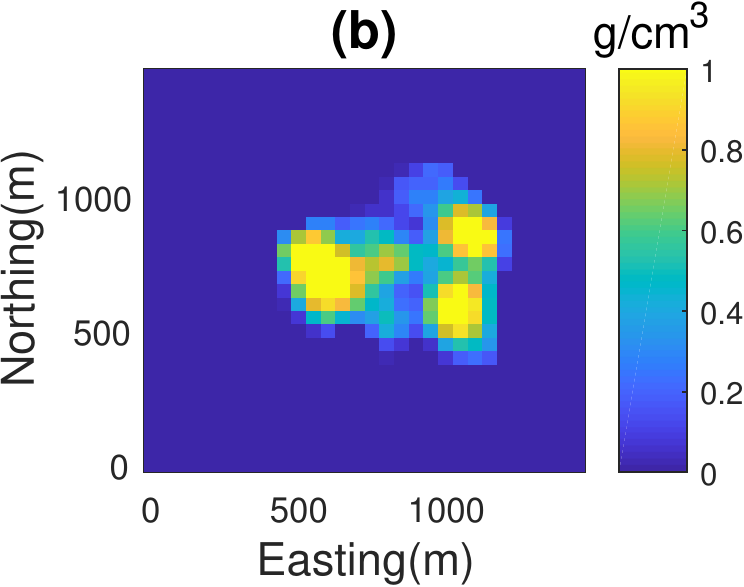}}
\caption {The reconstructed model for data in Fig.~\ref{fig2} using Algorithm~\ref{IterativeTVRGSVD} with $q=100$. (a) Cross-section at  $\mathrm{northing}=725$~m; (b) Plane-section at $\mathrm{depth}=100$~m.} \label{fig3}
\end{figure*}

\begin{figure*}
\subfigure{\label{fig4a}\includegraphics[width=.45\textwidth]{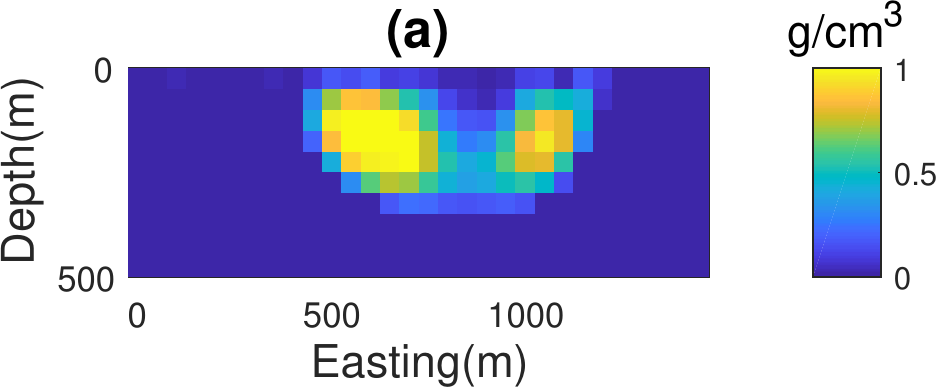}}
\subfigure{\label{fig4b}\includegraphics[width=.35\textwidth]{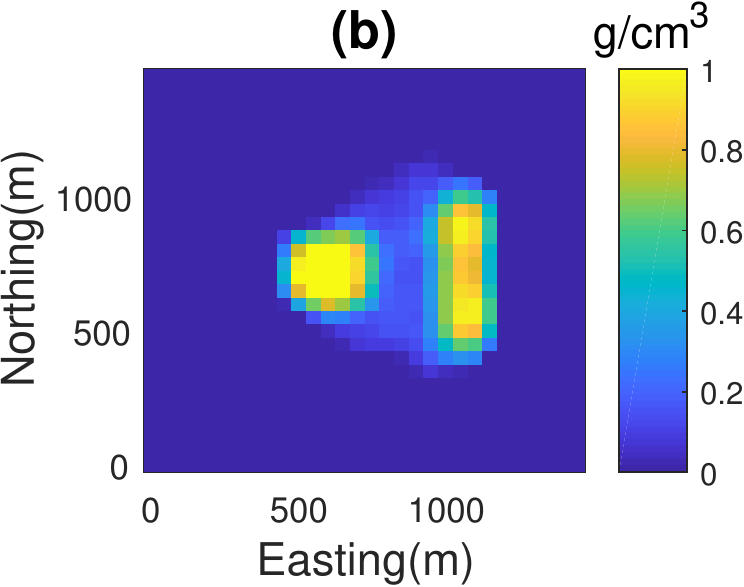}}
\caption {The reconstructed model for data in Fig.~\ref{fig2} using Algorithm~\ref{IterativeTVRGSVD} with $q=300$. (a) Cross-section at  $\mathrm{northing}=725$~m; (b) Plane-section at $\mathrm{depth}=100$~m.} \label{fig4}
\end{figure*}

\begin{figure*}
\subfigure{\label{fig5a}\includegraphics[width=.45\textwidth]{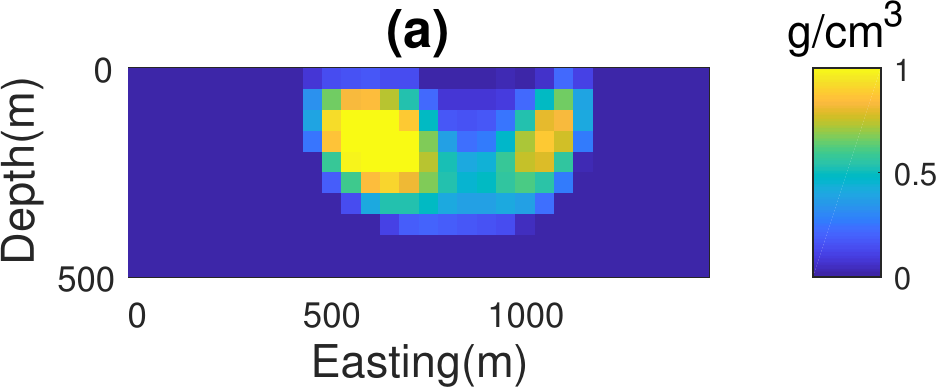}}
\subfigure{\label{fig5b}\includegraphics[width=.35\textwidth]{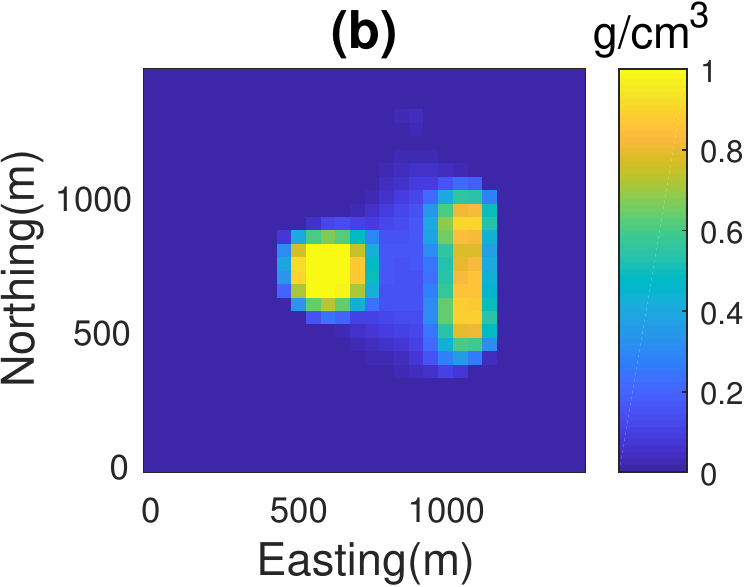}}
\caption {The reconstructed model for data in Fig.~\ref{fig2} using Algorithm~\ref{IterativeTVRGSVD} with $q=500$. (a) Cross-section at  $\mathrm{northing}=725$~m; (b) Plane-section at $\mathrm{depth}=100$~m.} \label{fig5}
\end{figure*}

\begin{figure*}
\subfigure{\label{fig6a}\includegraphics[width=.45\textwidth]{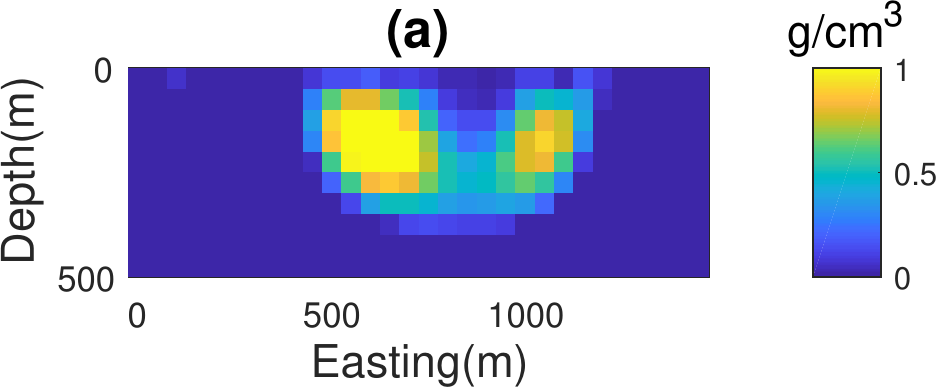}}
\subfigure{\label{fig6b}\includegraphics[width=.35\textwidth]{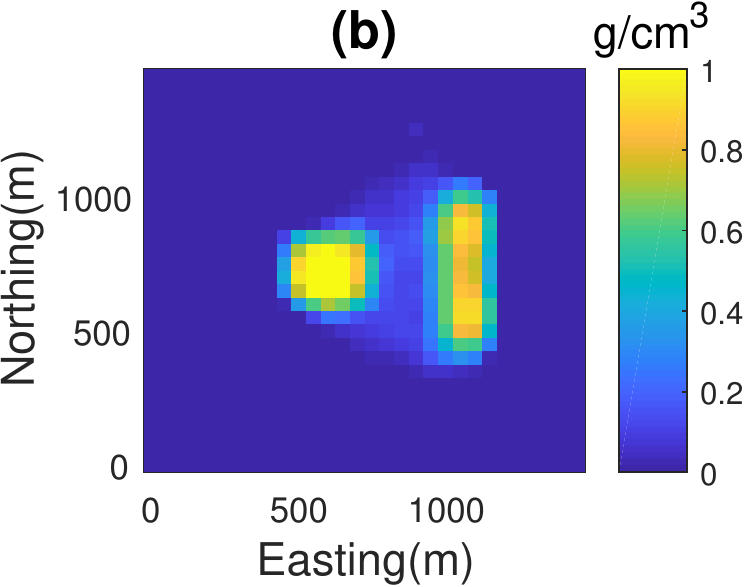}}
\caption {The reconstructed model for data in Fig.~\ref{fig2} using Algorithm~\ref{IterativeTVRGSVD} with $q=300$ at $K=50$. (a) Cross-section at  $\mathrm{northing}=725$~m; (b) Plane-section at $\mathrm{depth}=100$~m.} \label{fig6}
\end{figure*}

\begin{table}
\caption{Parameters and results for Algorithm~\ref{IterativeTVRGSVD} applied to the small scale dipping dikes model. The number of iterations is $K$, $\alpha^{(K)}$ is the final regularization parameter, $RE^{(K)}=\frac{\|\bfm_{\mathrm{exact}}-\bfm ^{(K)} \|_2}{\|\bfm_{\mathrm{exact}} \|_2}$ is the relative error of the reconstructed model at iteration $K$. The final $\chi_{\text{computed}}^2$ is  also reported.}\label{tab1}
\begin{tabular}{c  c  c  c  c | c c c cc}
\hline
\multicolumn{5}{c}{Input Parameters}&\multicolumn{5}{c}{Results}\\ \hline
$\bfma$&$ \rho_{\mathrm{min}}$ (g~cm$^{-3})$ &$\rho_{\mathrm{max}}$ (g~cm$^{-3}$)&$K_{\mathrm{max}}$&$q$&  $RE^{(K)}$&$\alpha^{(K)}$&  $K$ & $\chi_{\text{computed}}^2$  & Time (s)  \\ \hline
$\mathbf{0}$& $0$&$1$&$200$\\
&&&&$100$ & $0.7745$& $26.63$& $200$& $3978.9$& $494$ \\ 
&&&&$300$ & $0.7143$& $41.62$& $174$& $942.4$& $595$\\ 
&&&&$500$ & $0.7244$& $51.62$& $49$& $942.2$& $249$\\   \hline
$\mathbf{0}$& $0$&$1$&$50$ &$300$ & $0.7308$& $38.01$& $50$& $1165.7$& $165$\\  \hline
$\mathbf{0}$& $0$&$2$&$200$&$500$ & $0.7257$& $31.60$& $49$& $941.6$& $227$\\  \hline
$\neq \mathbf{0}$& $0$&$1$&$200$&$500$ & $0.7076$& $76.74$& $200$& $1952.6$& $972$\\  \hline
\end{tabular}
\end{table}
The parameters of the inversion needed for Algorithm~\ref{IterativeTVRGSVD}, and the results, are detailed for each example in Table~\ref{tab1}. These are the maximum number of iterations for the inversion  $K_{\mathrm{max}}$, the bound constraints  for the model  $\rho_{\mathrm{min}}$ and $\rho_{\mathrm{max}}$, the initial data $\bfma$ and the choice for $q$. Further the inversion terminates if $\chi_{\text{computed}}^2 =\| \Wd(\bfdo-\bfdp^{(K)})\|_2^2  \leq m+\sqrt{2m}=942.4$ is satisfied for iteration $K<K_{\mathrm{max}}$. The results for different choices of $q$, but all other parameters the same, are illustrated in Figs.~\ref{fig3}, \ref{fig4} and  \ref{fig5}, respectively.  With $q=100$ two dipping structures are recovered but  the iteration has not converged by $K=200$, and both the relative error and   $\chi_{\text{computed}}^2$ are large. With $q=300$ and $500$ better reconstructions of the dikes are achieved although the extension of the left dike is overestimated. While the errors are nearly the same,  the inversions terminate at $K=174$ and $K=49$ for $q=300$ and $q=500$, respectively. This leads to low computational time when using $q=500$ as compared with the other two cases, see Table~\ref{tab1}. We should note here that for all three cases the relative error decreases rapidly for the early iterations, after which there is little change in the model between iterations. For example the result for $q=300$ at iteration $K=50$, as illustrated in  Fig.~\ref{fig6} and detailed in Table~\ref{tab1}, is acceptable and is achieved with a substantial reduction in the computational time.

\begin{figure*}
\subfigure{\label{fig7a}\includegraphics[width=.45\textwidth]{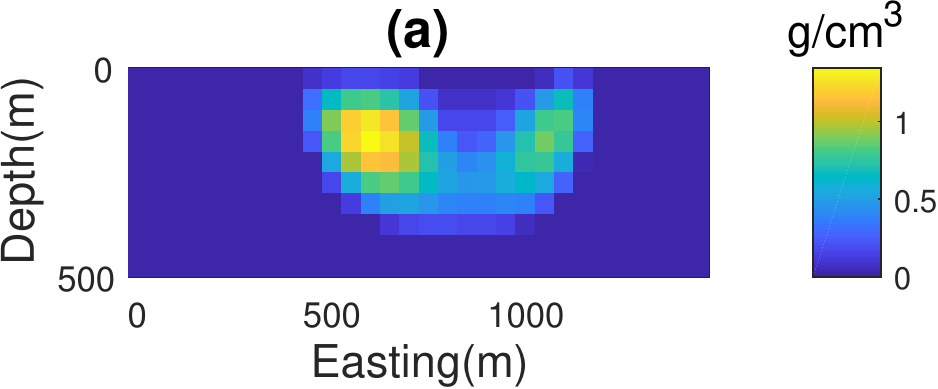}}
\subfigure{\label{fig7b}\includegraphics[width=.35\textwidth]{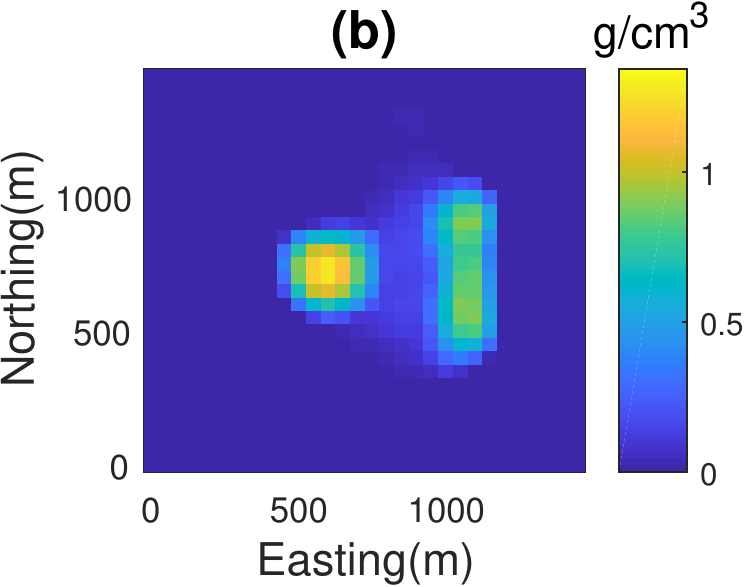}}
\caption {The reconstructed model for data in Fig.~\ref{fig2} using Algorithm~\ref{IterativeTVRGSVD} with $q=500$ when uncorrected upper density bound $\rho_{\mathrm{max}}=2$~g~cm$^{-3}$ is selected. (a) Cross-section at  $\mathrm{northing}=725$~m; (b) Plane-section at $\mathrm{depth}=100$~m.} \label{fig7}
\end{figure*}

\begin{figure*}
\subfigure{\label{fig8a}\includegraphics[width=.45\textwidth]{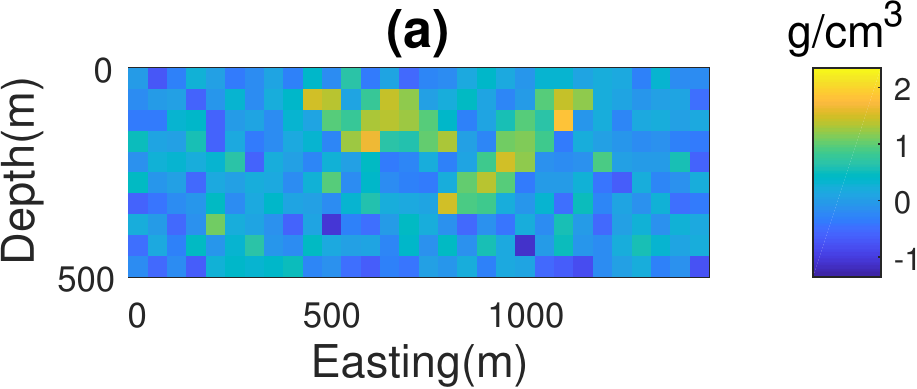}}
\subfigure{\label{fig8b}\includegraphics[width=.45\textwidth]{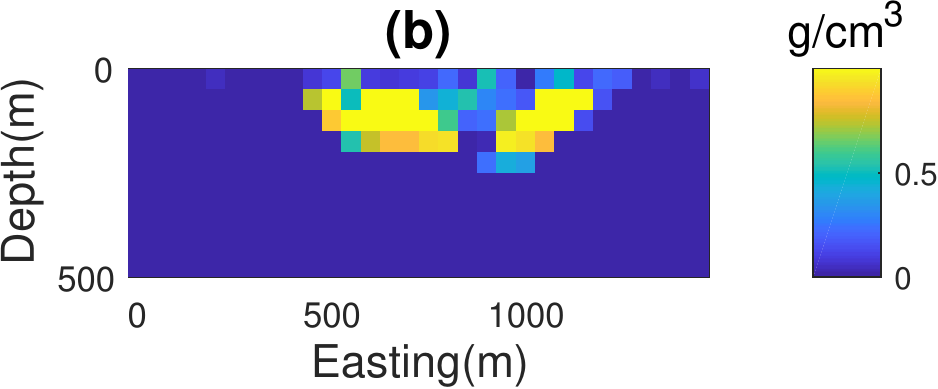}}
\caption { (a) The initial model which generated by adding Gaussian noise with standard deviation of $ (0.05~\bfm_{\mathrm{true}} + 0.02~\| \bfm_{\mathrm{true}} \|)$ to the true model. (b) The reconstructed model for data in Fig.~\ref{fig2} using Algorithm~\ref{IterativeTVRGSVD} with $q=500$ when model shown in Fig.~\ref{fig8a} is used as $\bfma$.} \label{fig8}
\end{figure*}

As compared to inversion using the $L_2$-norm of the gradient of the model parameters, see for example Li \& Oldenburg \shortcite{LiOl:96}, or conventional minimum structure inversion, the results obtained using Algorithm~\ref{IterativeTVRGSVD} provide a subsurface model that is not smooth. Further, as compared with $L_{0}$ and $L_{1}$ norms  applied for the model parameters, see for example \cite{LaKu:83,PoZh:99,VRA:2017a}, the TV inversion does not yield a model that is sparse or compact. On the other hand, the TV inversion is far less dependent on correct specification of the model constraints. This is illustrated in Fig.~\ref{fig7}, which is the same case as Fig.~\ref{fig5}, but with $\rho_{\mathrm{max}}=2$~g~cm$^{-3}$, and demonstrates that the approach is generally robust. Finally, to contrast with the algorithm presented in Bertete-Aguirre et al. \shortcite{BCO:2002} the results in Fig~\ref{fig8b} are for an alternative choice for $\bfma$, illustrated in Fig.~\ref{fig8a}, and $q=500$ which is consistent with the approach in Bertete-Aguirre et al. \shortcite{BCO:2002}. Here  $\bfma$ is obtained by taking the true model with the addition of Gaussian noise with a standard deviation of $ (0.05~\bfm_{\mathrm{true}} + 0.02~\| \bfm_{\mathrm{true}} \|)$. Not surprisingly the reconstructed density model is more focused  and is closer to the true model. The computed value for  $\chi_{\text{computed}}^2$ is, however, larger than the specified value and the algorithm terminates at $K_{\mathrm{max}}$. This occurs due to the appearance of an incorrect density distribution in the first layer of the subsurface. Together, these results demonstrate that the TVRGSVD technique is successful and offers a good option for the solution of larger scale models. 

\subsection{Model of multiple bodies}\label{multiplebodies}
\begin{figure*}
\subfigure{\label{fig9a}\includegraphics[width=.45\textwidth]{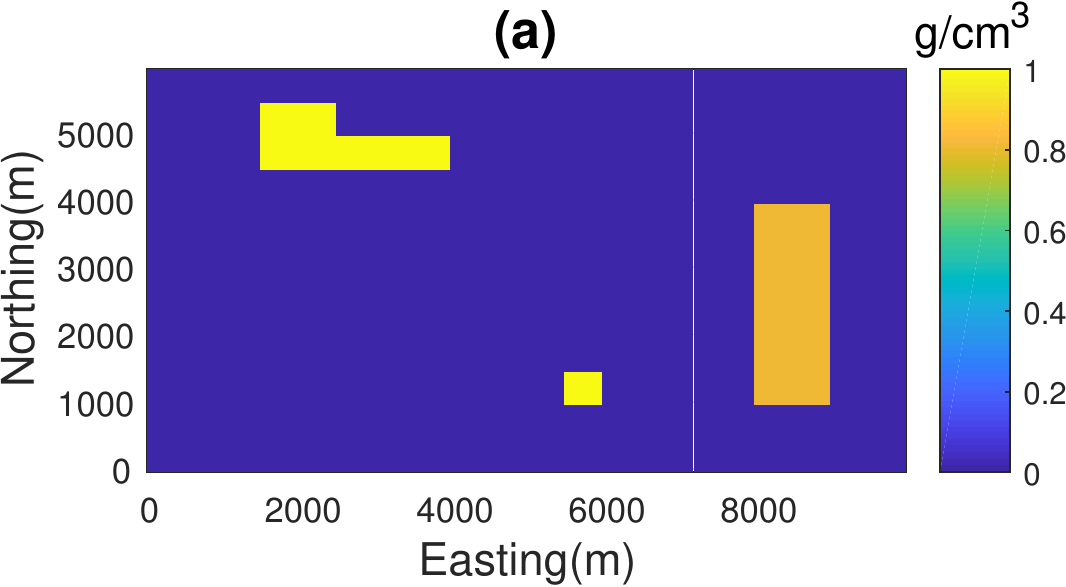}}
\subfigure{\label{fig9b}\includegraphics[width=.45\textwidth]{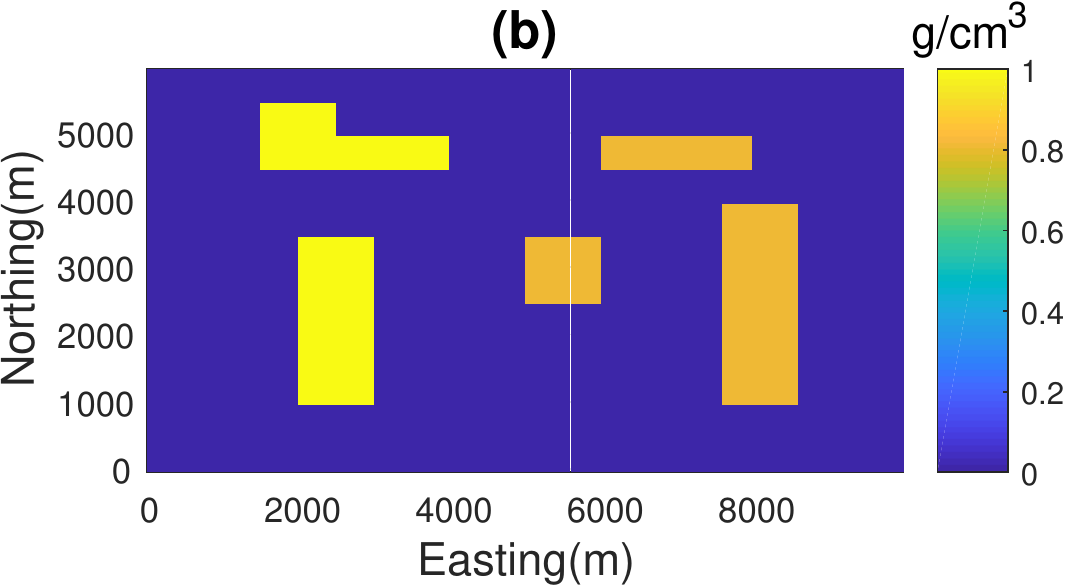}}
\subfigure{\label{fig9c}\includegraphics[width=.45\textwidth]{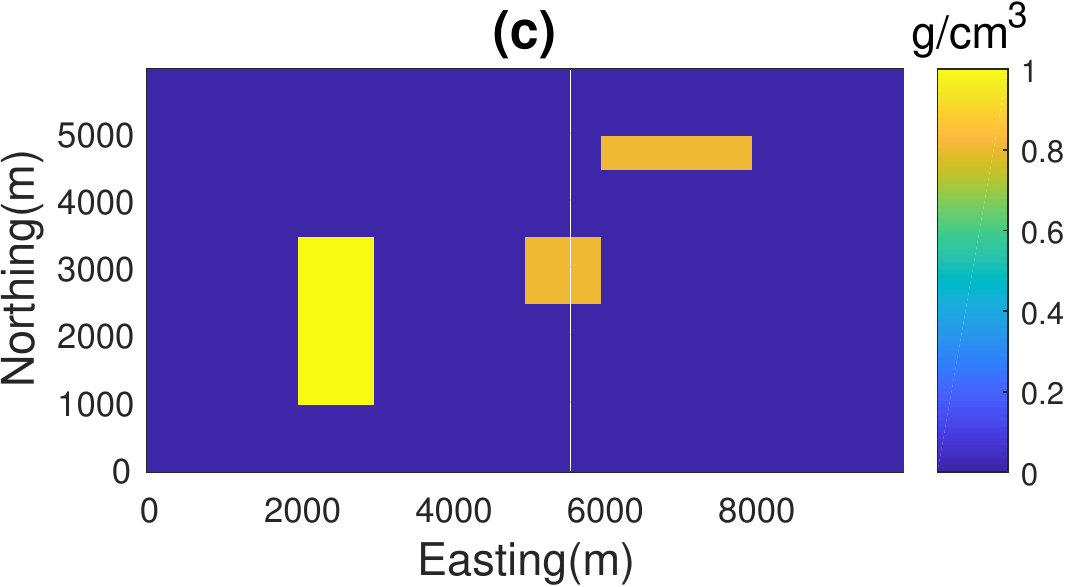}}
\subfigure{\label{fig9d}\includegraphics[width=.45\textwidth]{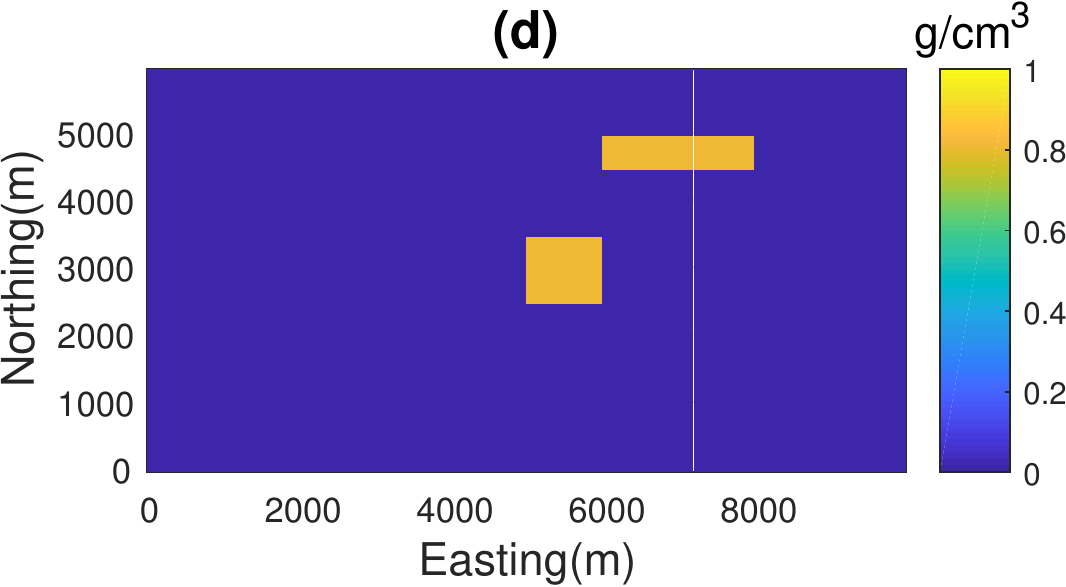}}
\caption { Model consists of six bodies with various geometries and sizes embedded in a homogeneous background. Bodies have the densities  $1$~g~cm$^{-3}$  and  $0.8$~g~cm$^{-3}$. (a) Plane-section at $\mathrm{depth}=100$~m; (b) Plane-section at $\mathrm{depth}=300$~m; (c) Plane-section at $\mathrm{depth}=500$~m; (d) Plane-section at $\mathrm{depth}=700$~m.} \label{fig9}
\end{figure*}

\begin{figure*}
\includegraphics[width=.5\textwidth]{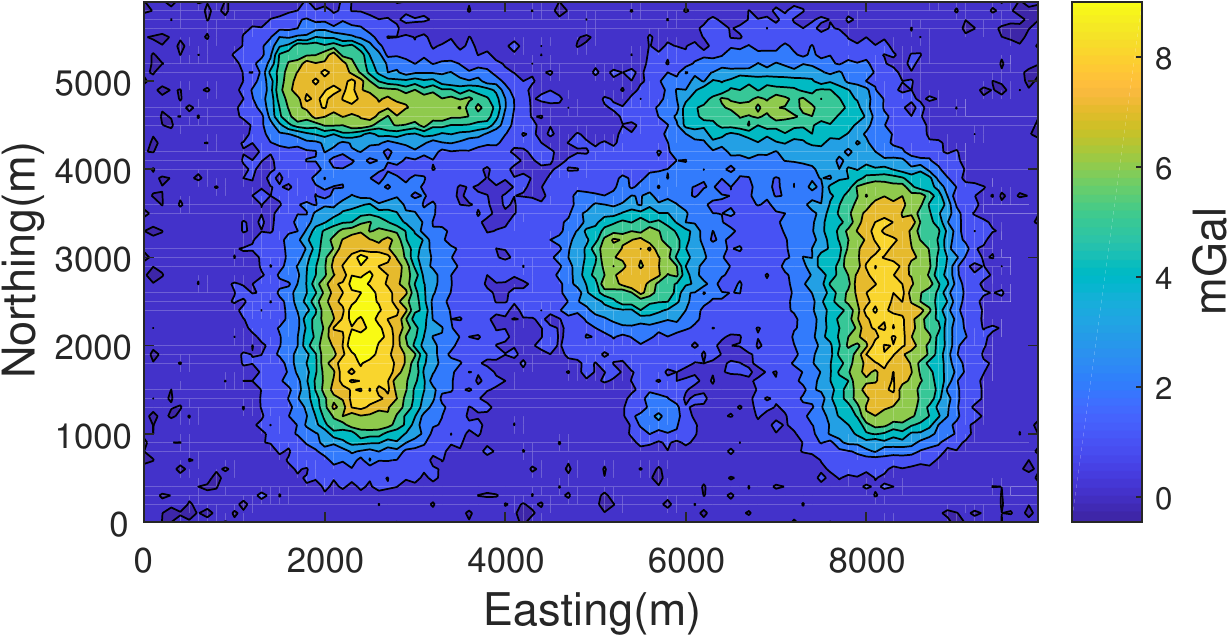}
\caption {The gravity anomaly produced by the multiple model shown in Fig.~\ref{fig9} and contaminated by Gaussian noise.} \label{fig10}
\end{figure*}

\begin{figure*}
\subfigure{\label{fig11a}\includegraphics[width=.45\textwidth]{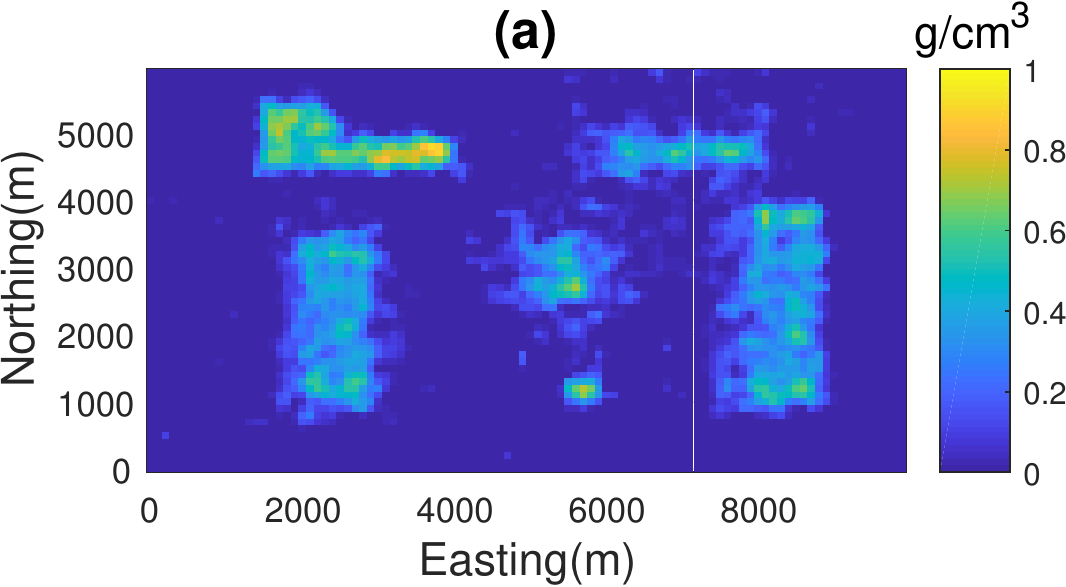}}
\subfigure{\label{fig11b}\includegraphics[width=.45\textwidth]{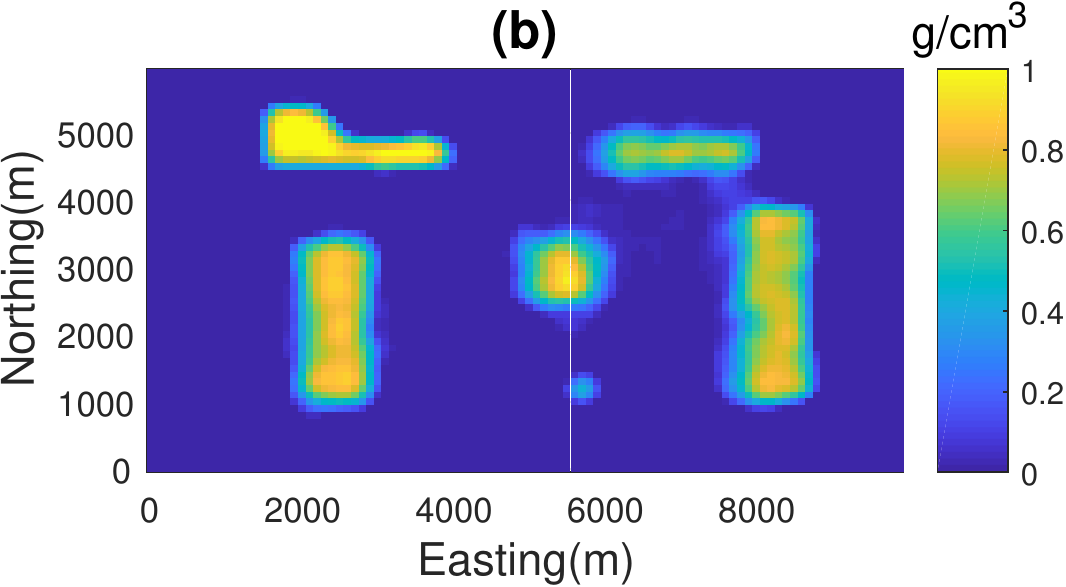}}
\subfigure{\label{fig11c}\includegraphics[width=.45\textwidth]{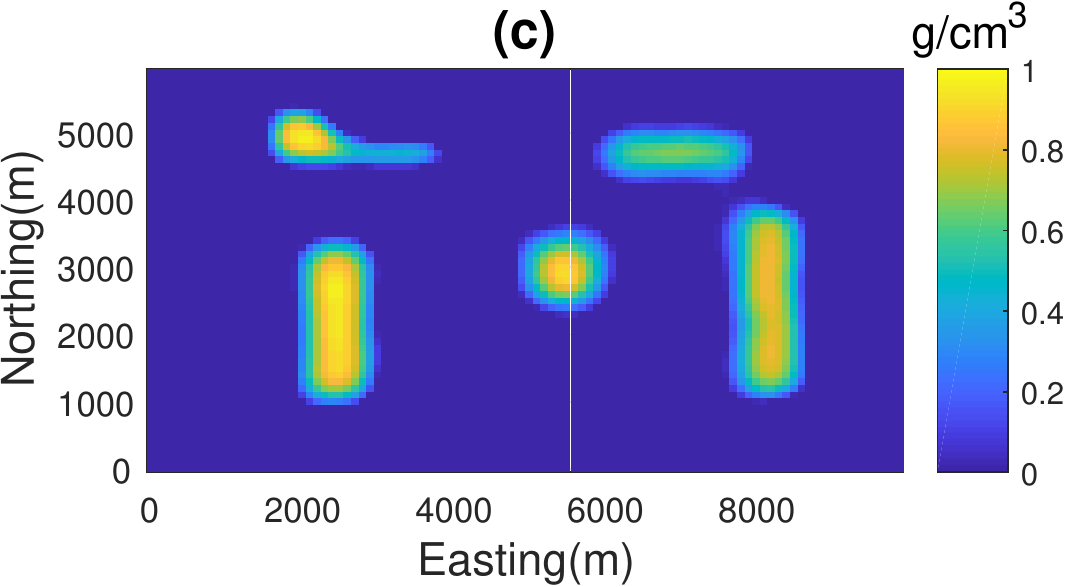}}
\subfigure{\label{fig11d}\includegraphics[width=.45\textwidth]{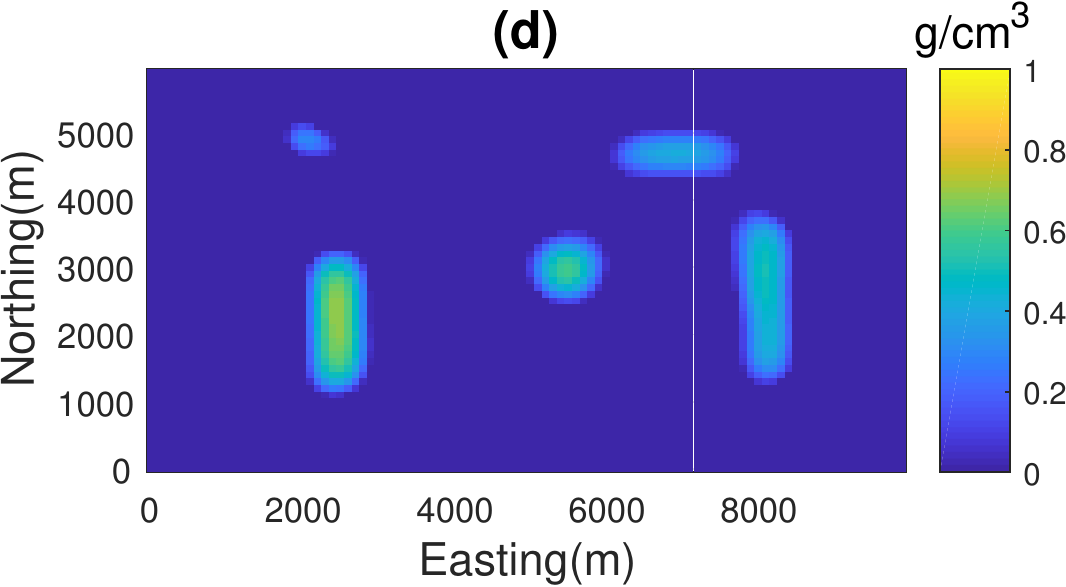}}
\caption { The reconstructed model for data in Fig.~\ref{fig10} using Algorithm~\ref{IterativeTVRGSVD} with $q=1000$. (a) Plane-section at $\mathrm{depth}=100$~m; (b) Plane-section at $\mathrm{depth}=300$~m; (c) Plane-section at $\mathrm{depth}=500$~m; (d) Plane-section at $\mathrm{depth}=700$~m.} \label{fig11}
\end{figure*}

\begin{figure*}
\includegraphics[width=.5\textwidth]{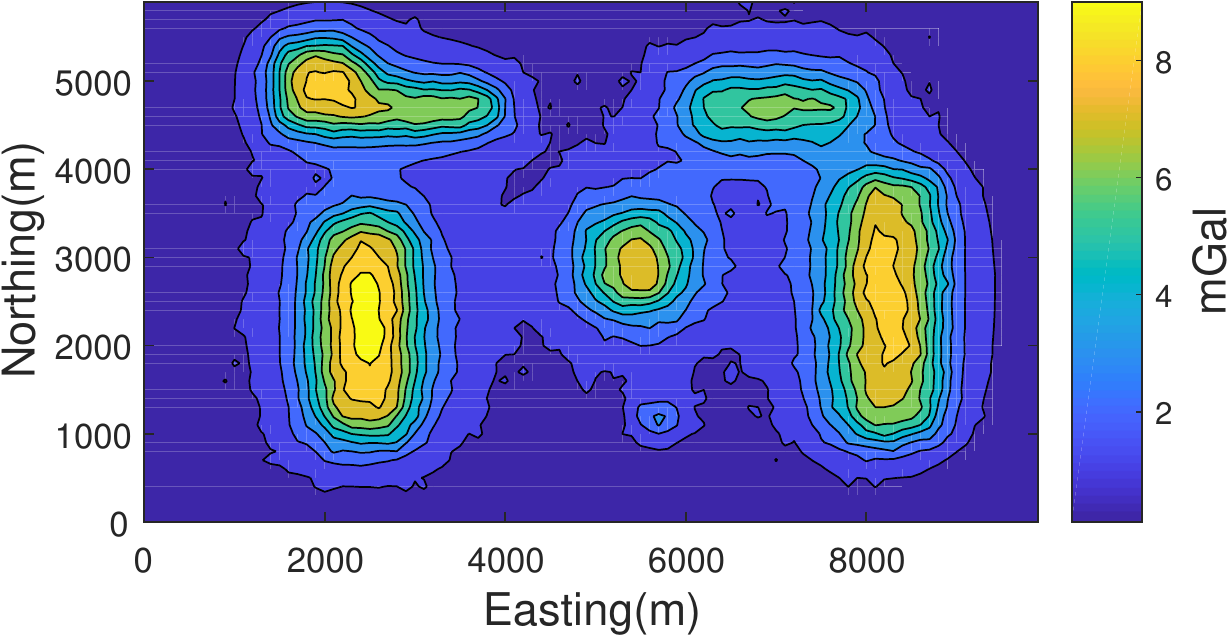}
\caption {The gravity anomaly produced by the reconstructed model shown in Fig.~\ref{fig11}.} \label{fig12}
\end{figure*}

We now consider a more complex and larger model that consists of six bodies with various geometries, sizes, depths and densities, as illustrated in Fig.~\ref{fig9}, and is further detailed in  Vatankhah et al \shortcite{VRA:2017a}. The surface gravity data are calculated on a $100 \times 60$ grid with $100$~m spacing. Gaussian noise with standard deviation of $ (0.02~(\bfde)_i + 0.001~\| \bfde \|)$ is added to each datum, giving the noisy data set as shown in Fig.~\ref{fig10}. The subsurface is divided into $100 \times 60 \times 10 = 60000$  cubes with sizes $100$~m in each dimension.  For the inversion we use Algorithm~\ref{IterativeTVRGSVD} with $\bfma=\mathbf{0}$, $\rho_{\mathrm{min}}=0$~g~cm$^{-3}$, $\rho_{\mathrm{max}}=1$~g~cm$^{-3}$ and $K_{\mathrm{max}}=50$. We perform the inversion for three values of $q$, $500$, $1000$ and  $2000$, and report the results in Table~\ref{tab2}. The reconstruction with $q=500$ is less satisfactory than that achieved with larger $q$ and the iterations terminate at $K_{\mathrm{max}}$ with a large value of $\chi_{\text{computed}}^2$.  The results with $q=1000$ and $q=2000$ have similar relative errors, but the computational cost is much reduced using $q=1000$; although the desired $\chi_{\text{computed}}^2$  is not achieved the result is close and acceptable. The reconstructed model and associated gravity response, using $q=1000$, are illustrated in Figs.~\ref{fig11} and~\ref{fig12},  respectively. While the maximum depths of the anomalies are overestimated, the horizontal borders are reconstructed accurately. These two examples suggest that  $q > m/6$ is suitable for Algorithm~\ref{IterativeTVRGSVD}, which confirms our previous conclusions when using the randomized SVD, see Vatankhah et al. \shortcite{VRA:2017b}. We should note that as compared to the case in which the standard Tikhonov functional can be used, see such results with the focusing inversion in Vatankhah et al. \shortcite{VRA:2017a,VRA:2017b}, the computational cost here is much higher.

\begin{table}
\caption{The results of the inversion for the multiple bodies example using Algorithm~\ref{IterativeTVRGSVD}.}\label{tab2}
\begin{tabular}{c  c  c  c  c | c c c cc}
\hline
\multicolumn{5}{c}{Input Parameters}&\multicolumn{5}{c}{Results}\\ \hline
$\bfma$&$ \rho_{\mathrm{min}}$ (g~cm$^{-3})$ &$\rho_{\mathrm{max}}$ (g~cm$^{-3}$)&$K_{\mathrm{max}}$&$q$&  $RE^{(K)}$&$\alpha^{(K)}$&  $K$ & $\chi_{\text{computed}}^2$  & Time (s)  \\ \hline
$\mathbf{0}$ &$0$&$1$&$50$&$500$ & $0.6710$& $15.31$& $50$& $11719.2$& $6014$\\   \hline
$\mathbf{0}$ &$0$&$1$&$50$&$1000$ & $0.6451$& $18.02$& $50$& $7324.0$& $7262$\\   \hline
$\mathbf{0}$ &$0$&$1$&$50$&$2000$ & $0.6452$& $18.79$& $44$& $6103.4$& $10666$\\   \hline
\end{tabular}
\end{table}

\section{Real data}\label{real}
\begin{figure*}
\subfigure{\label{fig13a}\includegraphics[width=.35\textwidth]{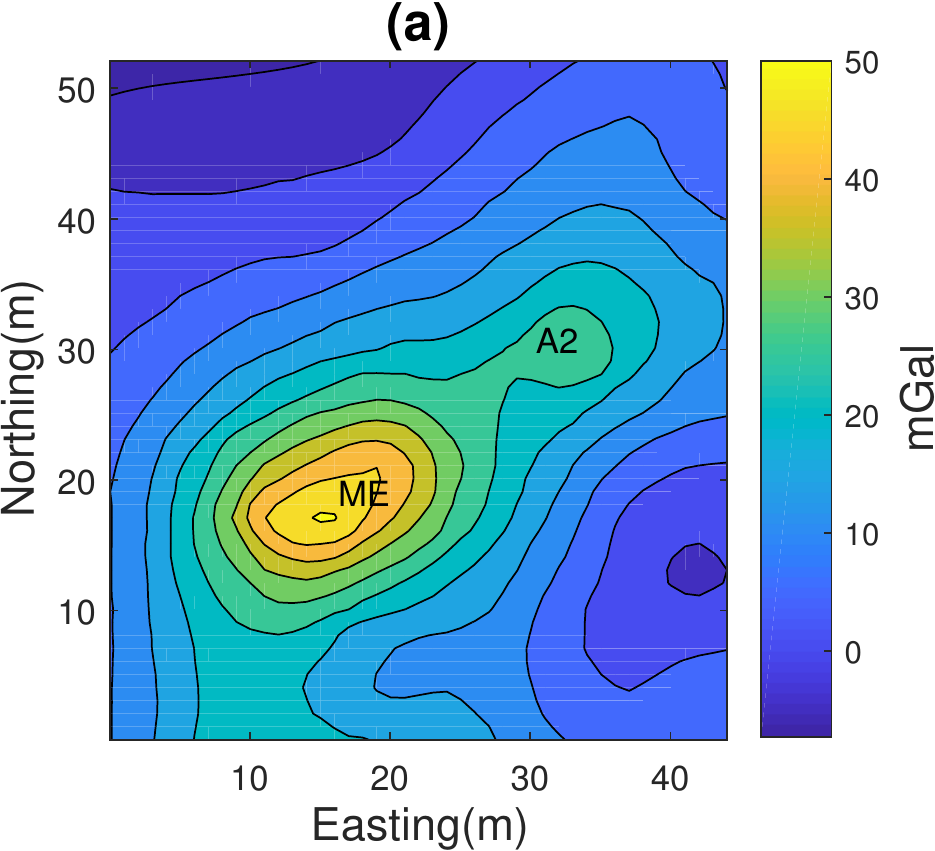}}
\subfigure{\label{fig13b}\includegraphics[width=.35\textwidth]{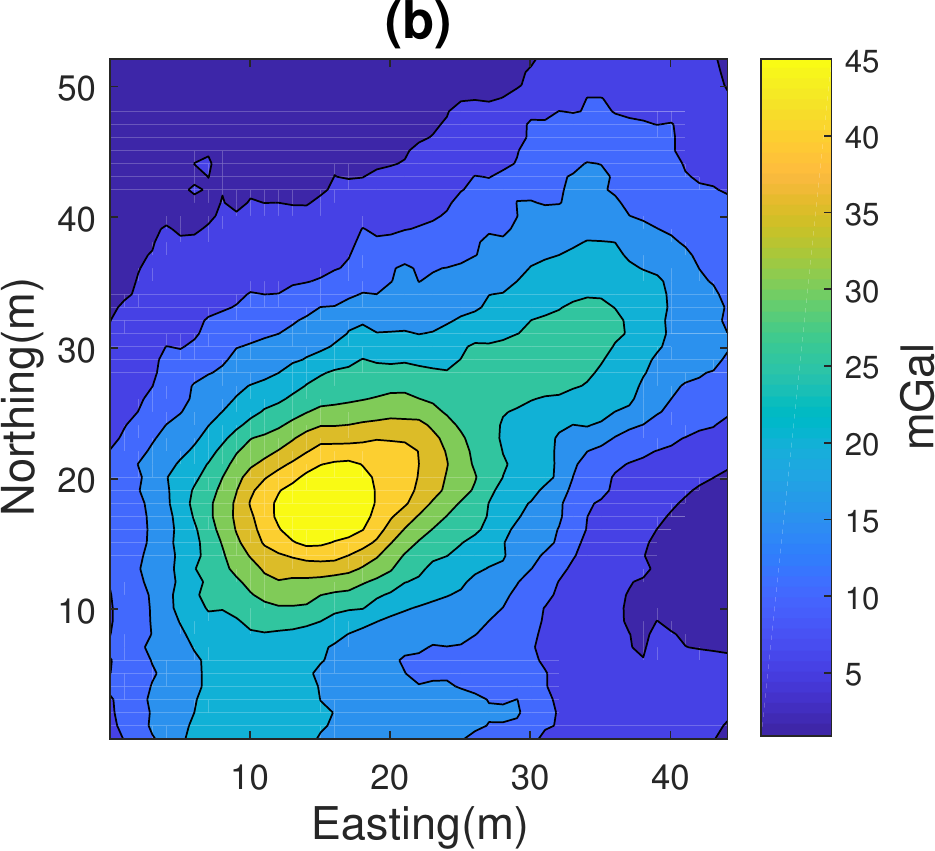}}
\caption {(a) Residual gravity data over the Morro do Engenho complex digitized from Dutra \& Marangoni \shortcite{DuMa:2009}; (b) The response of the reconstructed model shown in Fig.~\ref{fig14}.} \label{fig13}
\end{figure*}

We use the gravity data from the Goi\'as Alkaline Province (GAP)  of the central region of  Brazil. The GAP is a result of mafic-alkaline magmatism that occurred in the Late Cretaceous  and includes mafic-ultramafic alkaline complexes in the northern portion, subvolcanic  alkaline intrusions in the central region and volcanic products to the south with several dikes throughout the area \cite{Du:2012}. We select a region from the northern part of GAP in which Morro do Engenho Complex (ME) outcrops and another intrusive body, A$2$, is completely covered by Quaternary sediments \cite{DuMa:2009}. The data was digitized carefully from Fig.~$3$ in Dutra \& Marangoni \shortcite{DuMa:2009} and re-gridded into $45 \times 53= 2385$ data points with spacing $1$~km, Fig.~\ref{fig13a}. For these data, the result of the smooth inversion using Li \& Oldenburg \shortcite{LiOl:98} algorithm was presented in Dutra \& Marangoni \shortcite{DuMa:2009} and Dutra et al. \shortcite{Du:2012}. Furthermore, the result of focusing inversion based on $L_1$-norm stabilizer is available in Vatankhah et al. \shortcite{VRA:2017b}. The results using the TV inversion presented here can therefore be compared with the inversions using both aforementioned algorithms.

\begin{figure*}
\subfigure{\label{fig14a}\includegraphics[width=.3\textwidth]{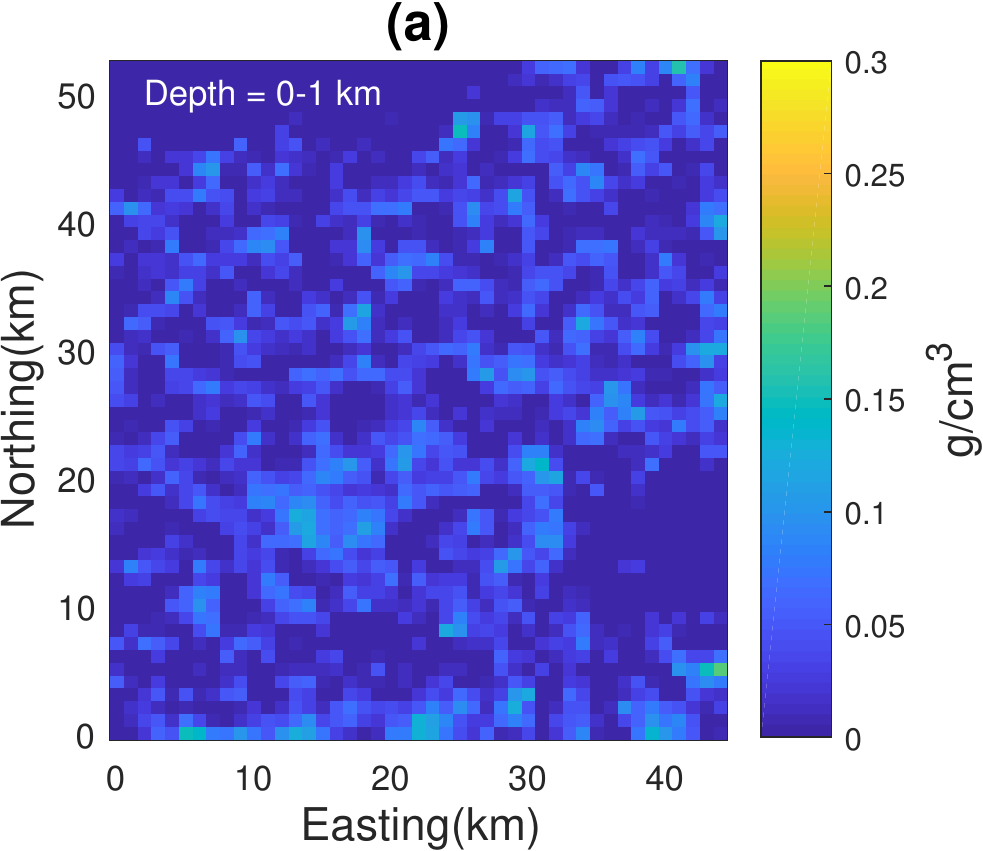}}
\subfigure{\label{fig14b}\includegraphics[width=.3\textwidth]{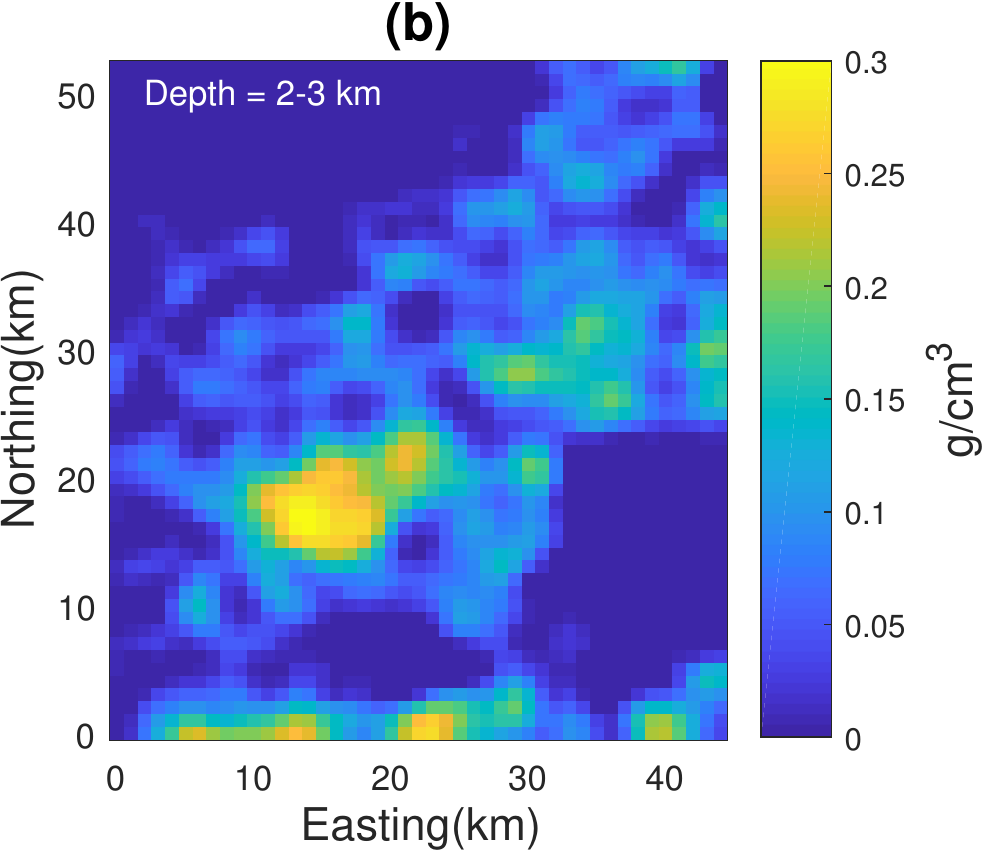}}
\subfigure{\label{fig14c}\includegraphics[width=.3\textwidth]{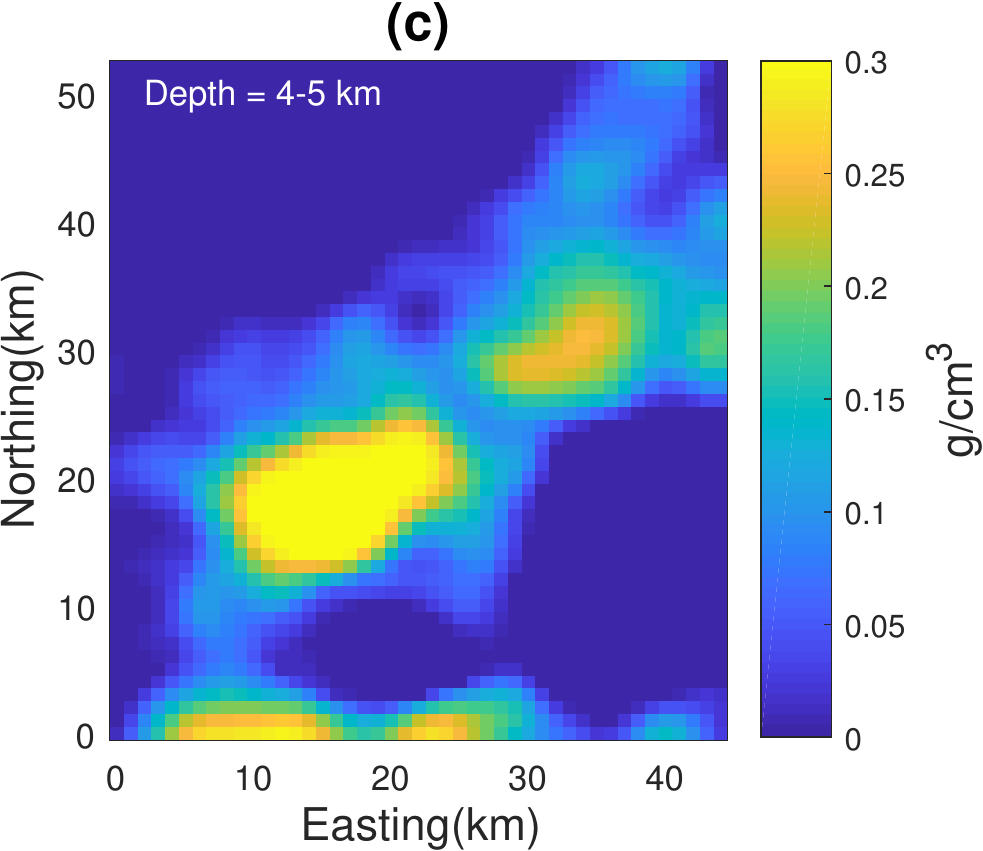}}
\subfigure{\label{fig14d}\includegraphics[width=.3\textwidth]{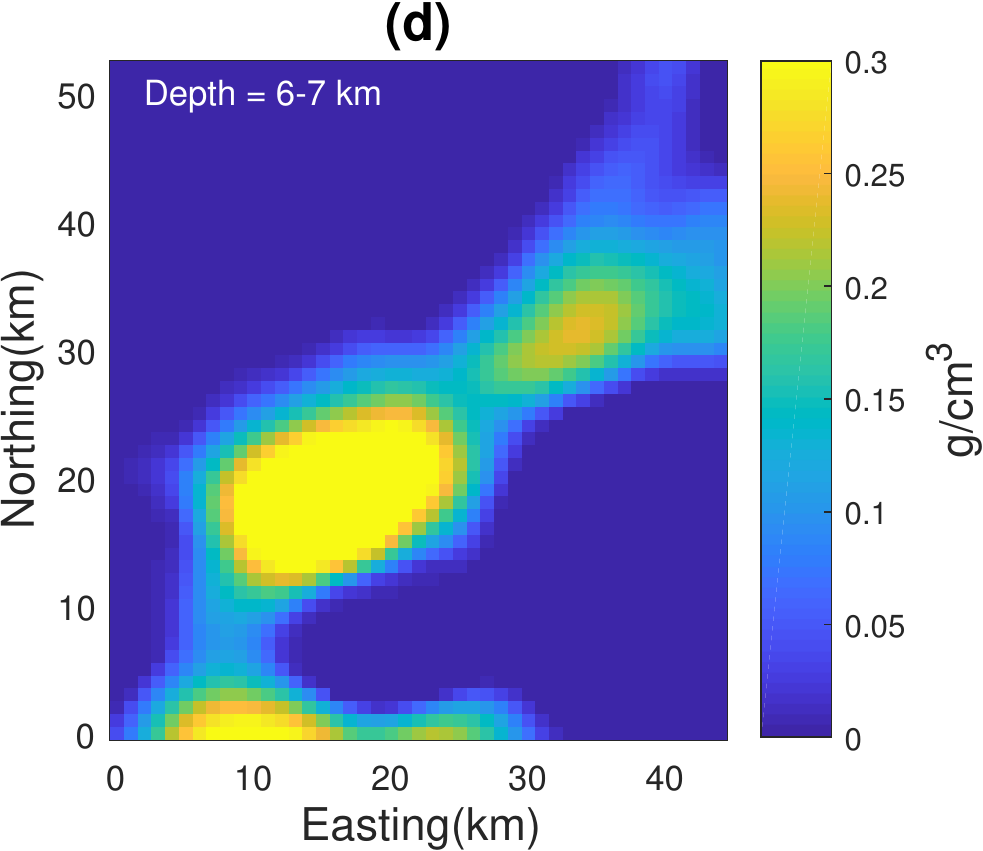}}
\subfigure{\label{fig14e}\includegraphics[width=.3\textwidth]{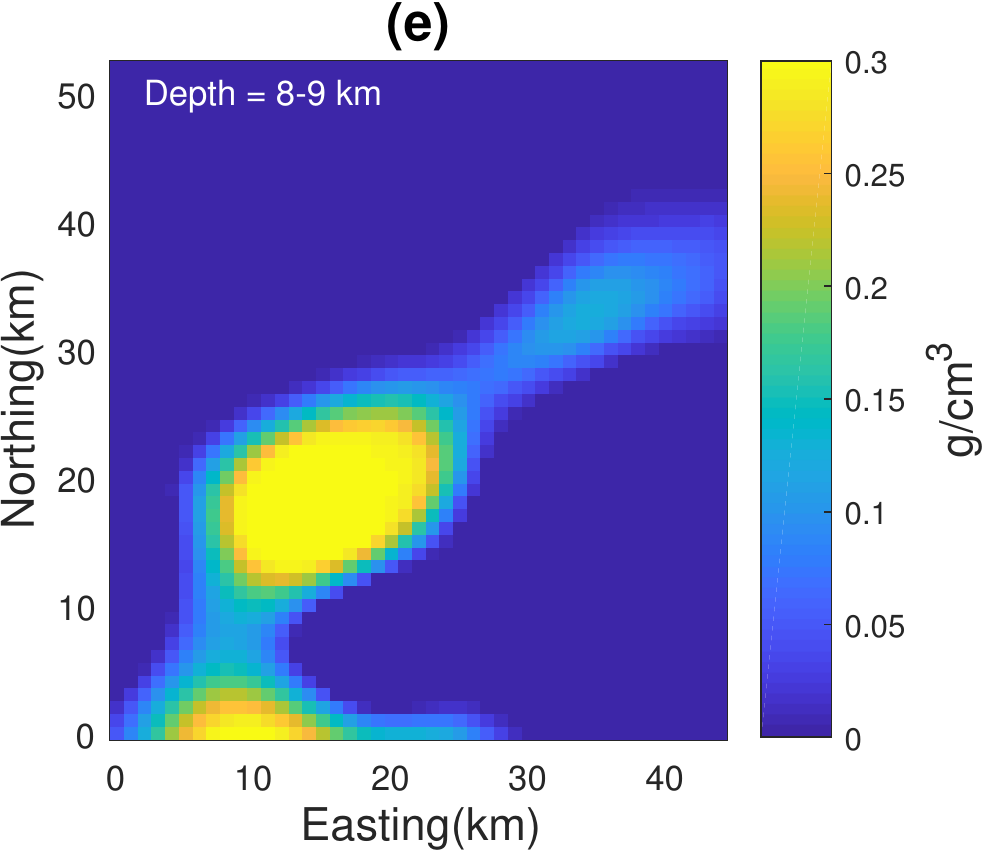}}
\subfigure{\label{fig14f}\includegraphics[width=.3\textwidth]{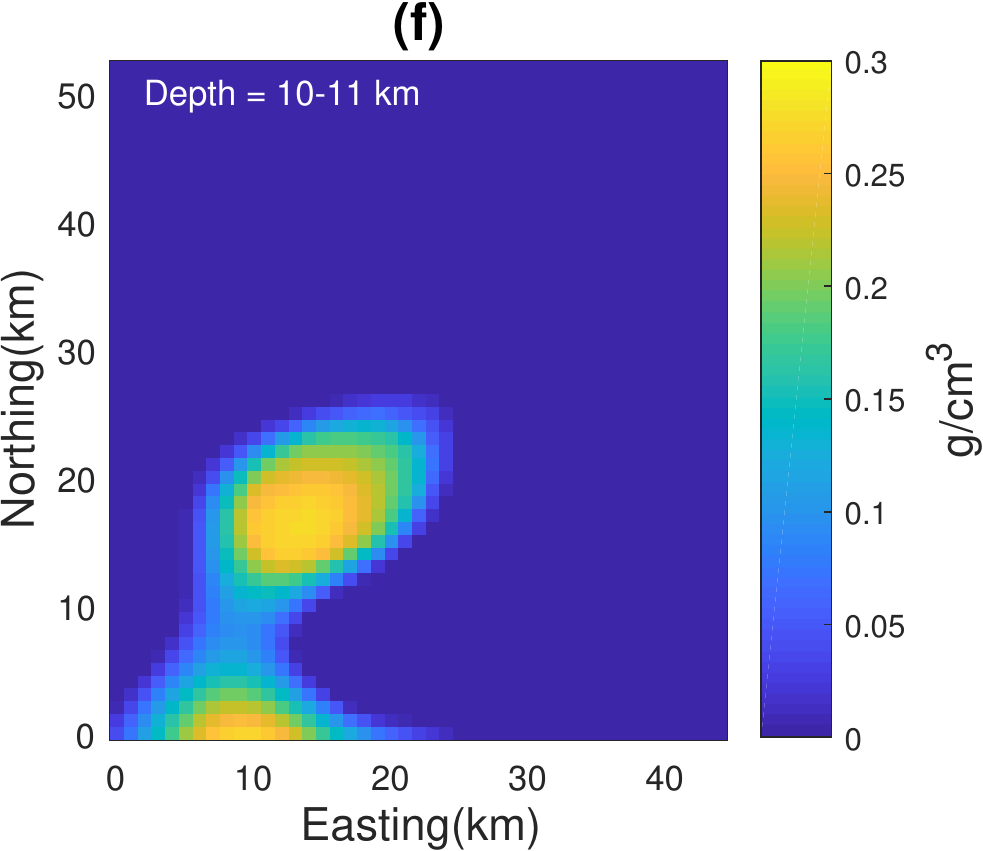}}
\caption { The plane-sections of the reconstructed model for the data in Fig.~\ref{fig13a} using
Algorithm~\ref{IterativeTVRGSVD} with q = 400. The sections are at the depths specified in the figures.} \label{fig14}
\end{figure*}

For the inversion we divide the subsurface into $45 \times 53 \times 14 =33390$ rectangular $1-$~km prisms. The density bounds  $\rho_{\mathrm{min}}=0$~g~cm$^{-3}$ and  $\rho_{\mathrm{max}}=0.3$~g~cm$^{-3}$  are selected based on geological information from Dutra \& Marangoni \shortcite{DuMa:2009}.  Algorithm~\ref{IterativeTVRGSVD} was implemented with $q=500$ and $K_{\mathrm{max}}=100$. The results of the inversion are presented in Table~\ref{tab3}, with the illustration of the predicted data due to the reconstructed model in Fig.~\ref{fig13b} and the reconstructed model in Fig.~\ref{fig14}.  As compared to the smooth estimate of the subsurface shown in  Dutra \& Marangoni \shortcite{DuMa:2009}, the obtained subsurface is blocky and non-smooth. The subsurface is also not  as sparse as that obtained by the focusing inversion in Vatankhah et al. \shortcite{VRA:2017b}. The ME and A$2$ bodies extend up to maximum $12$~km and $8$~km, respectively. Unlike the result obtained using the focusing inversion, for the TV inversion the connection between ME and A$2$ at depths $4$~km to $7$~km is not strong. We should note that the computational time for the focusing inversion presented in  Vatankhah et al. \shortcite{VRA:2017b} is much smaller than for the TV algorithm presented here.

\begin{table}
\caption{The results of the inversion for the data presented in Fig.~\ref{fig13a} using Algorithm~\ref{IterativeTVRGSVD}.}\label{tab3}
\begin{tabular}{c  c  c  c  c | c c c c}
\hline
\multicolumn{5}{c}{Input Parameters}&\multicolumn{4}{c}{Results}\\ \hline
$\bfma$&$ \rho_{\mathrm{min}}$ (g~cm$^{-3})$ &$\rho_{\mathrm{max}}$ (g~cm$^{-3}$)&$K_{\mathrm{max}}$&$q$&$\alpha^{(K)}$&  $K$ & $\chi_{\text{computed}}^2$  & Time (s)  \\ \hline
$\mathbf{0}$&$0$&$0.3$&$100$ &$500$& $102.48$& $38$& $2453.9$& $1533$\\   \hline
\end{tabular}
\end{table}

\section{Conclusions}\label{conclusion}
We developed an algorithm for total variation regularization applied to the  $3$-D gravity inverse problem. The presented algorithm provides a non-smooth and blocky image of the subsurface which may be useful when discontinuity of the subsurface is anticipated. Using the randomized generalized singular value decomposition, we have demonstrated that moderate-scale problems can be solved in a reasonable computational time. This seems to be the first use of the TV regularization using the RGSVD for gravity inversion. Our presented results use an Alternating Direction implementation to further improve the efficiency and reduce the memory demands of the problem.  The results show that there is less sensitivity to the provision of good bounds for the density values, and thus the TV may have a role to play for moderate scale problems where limited information on subsurface model parameters is available. It was of interest to investigate the $3$-D algorithm developed here in view of the results of Bertete-Aguirre et al. \shortcite{BCO:2002}  which advocated TV regularization in the context of $2$-D gravity inversion. We obtained results that are comparable with the simulations  presented in Bertete-Aguirre et al. \shortcite{BCO:2002} but here for much larger problems. In our simulations we have found that the computational time is much larger than that required for the focusing algorithm presented in Vatankhah et al. \shortcite{VRA:2017a,VRA:2017b}. Because it is not possible to transform the TV regularization to standard Tikhonov form, the much more expensive GSVD algorithm has to be used in place of the SVD that can be utilized for the focusing inversion presented in  Vatankhah et al. \shortcite{VRA:2017b}. On the other hand, the advantage of the TV inversion as compared to the focusing inversion is the lesser dependence on the density constraints and the generation of a subsurface which is not as smooth and admits discontinuities. The impact of the algorithm was illustrated for the inversion of real gravity data from the Morro do Engenho complex in central Brazil.

\begin{acknowledgments}
 
\end{acknowledgments}

\appendix
\section{The generalized singular value decomposition}\label{gsvd}
Suppose $\tilde{\tilde{G}} \in \mathcal{R}^{m\times n}$, $\tilde{D} \in \mathcal{R}^{p \times n}$ and  $\mathcal{N}(\tilde{\tilde{G}}) \cap \mathcal{N}(\tilde{D}) = 0$, where $\mathcal{N}(\tilde{\tilde{G}})$ is the null space of matrix $\tilde{\tilde{G}}$. Then there exist orthogonal matrices $U \in \mathcal{R}^{m\times m}$, $V \in \mathcal{R}^{p \times p}$ and a nonsingular matrix $X \in \mathcal{R}^{n \times n}$ such that  $\tilde{\tilde{G}}=U \Lambda X^T$ and $\tilde{D}=V M X^T$ \cite{PaSa:1981,ABT:2013}. Here, $\Lambda \in \mathcal{R}^{m \times n}$ is zero except for  entries $0< \Lambda_{1,(n-m)+1} \le \dots \Lambda_{m,n} <1$, and $M$ is diagonal of size $p\times n$ with entries $M_{1,1}> M_{2,2}\ge \dots \ge M_{p^*,p^*}>0$, where $p^*:=\mathrm{min}(p, n)$. The generalized singular values of the matrix pair $[\tilde{\tilde{G}}, \tilde{D}]$ are $\gamma_i=\lambda_i/\mu_i$, where $\gamma_1=\dots =\gamma_{(n-m)}=0<\gamma_{(n-m)+1}\le \dots \le \gamma_n$, and $\Lambda^T\Lambda =\mathrm{diag}(0,\dots 0,\lambda^2_{(n-m)+1},\dots, \lambda_n^2)$, $M^T M= \mathrm{diag}(1,\dots,1,\mu_{(n-m)+1}^2,\dots,\mu_n^2)$, and $\lambda_i^2+\mu_i^2=1$, $\forall i=1:n$, i.e. $M^T  M + \Lambda^T \Lambda=I_n$. 

Using the GSVD,  introducing $\bfu_{i}$ as the $i$th column of matrix $U$, we may immediately write the solution of \eqref{hsolution} as 
\begin{eqnarray}\label{gsvdsoln}
\bfh(\alpha)= \sum_{i=(n-m)+1}^{n} \frac{\gamma^2_i}{\gamma^2_i+\alpha^2} \frac{\bfu^T_{i-(n-m)}\tilde{\bfr}}{\lambda_i} (X^T)^{-1}_{i},
\end{eqnarray}
where $(X^T)^{-1}_{i}$ is the $i$th column of the inverse of the matrix $X^T$.

\end{document}